\documentclass[12pt,a4paper]{article}

\usepackage[english]{babel}

\usepackage{stmaryrd}
\usepackage{graphicx}
\usepackage[colorlinks=true, allcolors=blue]{hyperref}
\usepackage[english]{babel}
\usepackage{amsmath,amssymb,amsthm}
\usepackage{graphicx}
\usepackage{xspace}
\usepackage[colorlinks=true, allcolors=blue]{hyperref}
\usepackage{csquotes}

\usepackage{fancyhdr}          

\newcommand{\shorttitle}{Strong relations between $n$-surfaces, $n$-PCM's, and pseudomanifolds}

\newcommand{\myproof}{ \textbf{Proof: \ }}

\newcommand{\AC}{\alpha^{\square}\xspace}
\newcommand{\BC}{\beta^{\square}\xspace}
\newcommand{\TC}{\theta^{\square}\xspace}
\newcommand{\MAPPING}{\mathcal{F}\xspace}
\newcommand{\LINK}{\mathrm{lk}\xspace}
\newcommand{\DIM}{\mathrm{dim}\xspace}
\newcommand{\RANK}{\mathrm{rk}\xspace}
\newcommand{\PILOCAL}{\pi_{\mathrm{local}}\xspace}

\newcommand{\Reals}{\mathbb{R}\xspace}

\newcommand{\COMPOS}{\mathcal{CC}\xspace}
\newcommand{\INT}[2]{\llbracket #1,#2 \rrbracket\xspace}

\newtheorem{Theorem}{Theorem}
\newtheorem{Definition}{Definition}
\newtheorem{Proposition}{Proposition}
\newtheorem{Lemma}{Lemma}
\newtheorem{Notation}{Notation}

\title{Strong relations between discrete surfaces, poset-based connected manifolds, and normal pseudomanifolds}

\author{Nicolas Boutry\\
EPITA Research Laboratory (LRE)\\
14--16 rue Voltaire, FR-94270 Le Kremlin-Bicêtre, France\\
nicolas.boutry@lrde.epita.fr}

\begin{document}
\maketitle

\markboth{\shorttitle}{\shorttitle}  

\begin{abstract}
In this paper, we study some relationships existing between some particular mathematical structures: discrete surfaces coming from discrete topology and mathematical morphology, poset-based connected manifolds coming from discrete topology, and normal pseudomanifolds which are much used in discrete geometry and topological data analysis. We will also show that, even when poset-based connected manifolds are assumed to be simplicial complexes, and then supplied with many additional topological properties, they are not necessarily smooth. A set of sufficient conditions to ensure that poset-based connected manifolds are smooth will be provided.

\end{abstract}


\newpage

\section{Introduction}
\label{sec.intro}

In continuous topology~\cite{kelley2017general}, \textquote{pinches} in topological spaces can lead to issues (for example, they are not easy to parameterize due to the fact that they are not locally homeomorphic to $\Reals^n$).

\medskip

\begin{figure}[h!]
    \centering
    \includegraphics[width=0.25\linewidth]{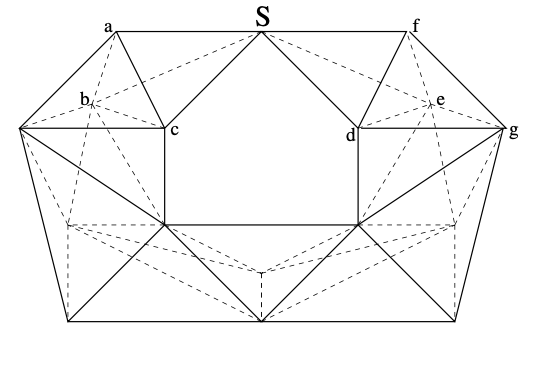}
    \caption{A pinched sphere where we see two $2$-faces that are sharing a vertex $S$ and thus are vertex-connected, but they are not locally edge-connected. From a discrete topology point of view, this simplicial complex is a pseudomanifold but not a normal pseudomanifold. At the same time, from a combinatorial/topological point of view, it is not a discrete surface. Additionally, from a continuous topology point of view, the geometric realization of this simplicial complex is not a topological manifold. The common reason of these three previous assertions is the pinch located at the vertex $S$. This figure is extracted from~\cite{daragon2005surfaces}.}
    \label{fig:pinchedsphere}
\end{figure}

In topology, pinches can be source of \emph{topological issues} (see Figure~\ref{fig:pinchedsphere}) in the sense that they lead to ambiguities in matter of local connectivities~\cite{latecki1995well,boutry.16.phd} in discrete settings (they can even lead to holes in surfaces using the marching cube algorithm~\cite{lorensen1987marching}), or in the sense that they make the structure they lie in difficult to parameterize~\cite{latecki19973d} in continuous settings.

\medskip

\begin{figure}[h!]
    \centering
    \includegraphics[width=0.25\linewidth]{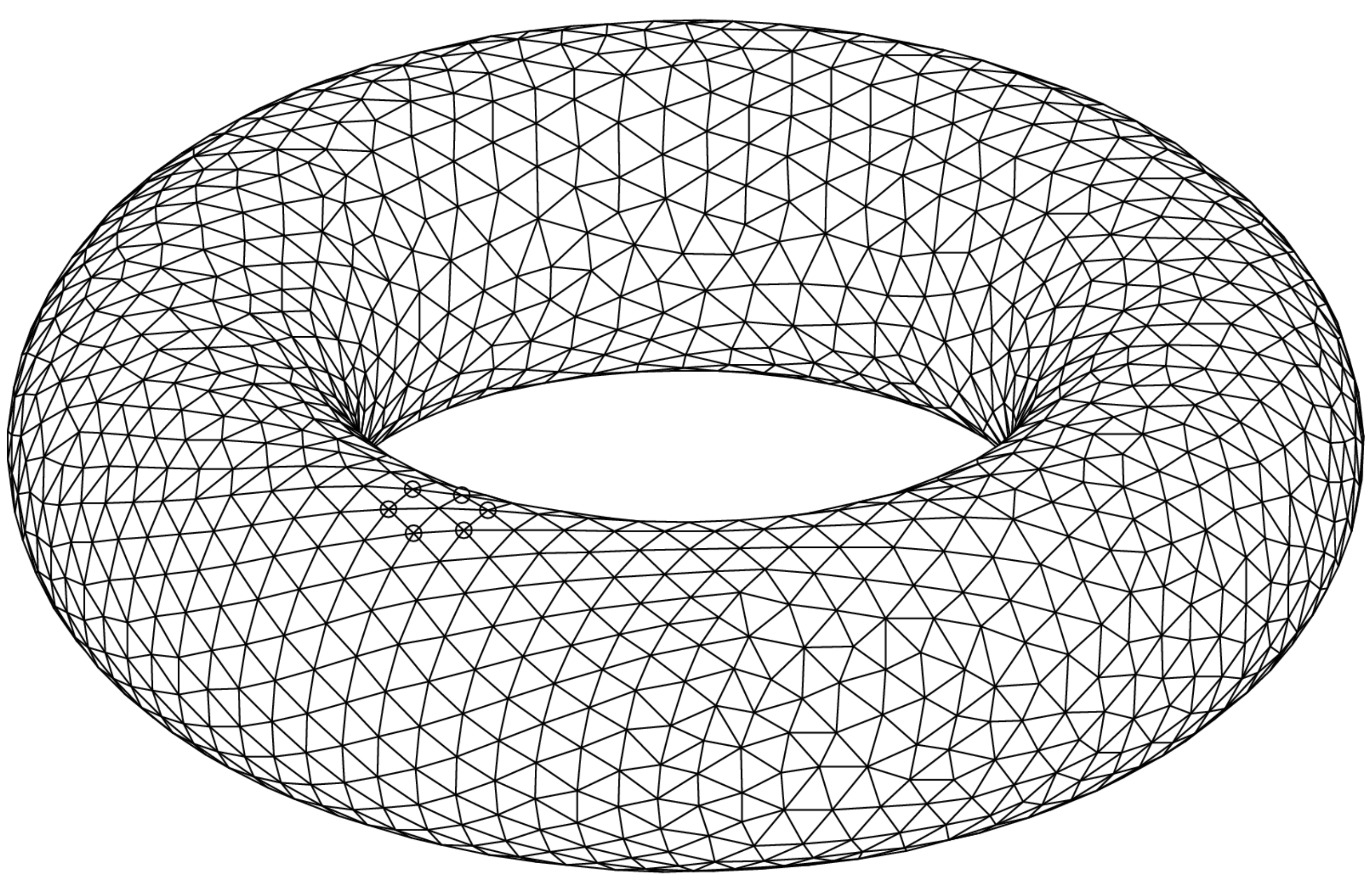}
    \caption{A topological structure which is at the same time a simplicial complex, a normal pseudomanifold, a topological manifold, a combinatorial $2$-manifold and a discrete $2$-surface. This figure is extracted from~\cite{boutry.16.phd}.}
    \label{fig:manifold}
\end{figure}

These topological issues explain why regular structures have been introduced in discrete topology, such as discrete surfaces~\cite{evako1995topological,daragon2005surfaces,najman2013discrete} or in the same family $n$-dimensional poset-based connected manifolds~\cite{boutry2024introducing} (shortly $n$-PCM's), or even \emph{normal} pseudomanifolds in algebraic topology~\cite{goresky1980intersection} (the links of their faces of codimension $2$ or more are pseudomanifolds). An example of simplicial complex which is at the same time a discrete $2$-surface, a normal $2$-pseudomanifold, and whose (underlying polyhedron of the) geometric realization is a topological manifold~\cite{lee2010introduction,latecki19973d,latecki1995well} is depicted in Figure~\ref{fig:manifold}.

\medskip

Many other kinds of manifolds can be found in the literature~\cite{munkres1984elements,lee2018introduction,ishihara1974quaternion,whitney1936differentiable}. In particular, combinatorial/stellar manifolds~\cite{daragon2005surfaces} are known to be even more regular than discrete surfaces (links of their vertices are spheres or balls) but we cannot always determine from a computational point of view if the studied structure is a combinatorial manifold or not, contrary to discrete surfaces, $n$-PCM's or even normal pseudomanifolds.

\medskip

This led the author to two questions: assuming that we start from a simplicial complex, are normal pseudomanifolds related to discrete surfaces and/or to $n$-PCM's since all these structures have the common point not to have \textquote{pinches}? Also, assuming that some $n$-PCM is also a simplicial complex, and then supplied with additional topological properties, does this $n$-PCM's becomes a smooth $n$-PCM? We answer to these different questions in this paper.

\medskip

We propose then to recall the mathematical background needed in this context in Section~\ref{sec.background}, to introduce some new properties relative to the studied structures in Section~\ref{sec.newbackground}, to present our main results relating discrete surfaces, poset-based connected manifolds, and normal pseudomanifolds in Section~\ref{sec.equivalence}, to show in Section~\ref{sec.pcm.vs.smooth.pcms} that $n$-PCM's that are simplicial complexes are not necessarily smooth $n$-PCMs, to conclude in Section~\ref{sec.conclusion}, and to present possible future works in Section~\ref{sec.future}. 

\section{Mathematical background}
\label{sec.background}

We will start this section with recalls in matter of axiomatic digital topology, so we will be able to define posets, discrete surfaces and then $n$-dimensional poset-based connected manifolds. We will continue with some relations that exist between them using the join operator. Finally, we will recall what are pseudomanifolds (especially normal pseudomanifolds).

\medskip

All over the paper, we will use the notation, with $a,b$ two integral values s.t. $a \leq b$, of \emph{integral intervals} $\INT{a}{b} := [a,b] \cap \mathbb{Z}$.

\subsection{Axiomatic digital topology}

\newcommand{\RC}{R^{\square}\xspace}

The recalls below are mainly extracted from~\cite{daragon2005discrete,najman2013discrete,boutry2020equivalence}. For $A$ and $B$ two sets of arbitrary elements, $A \times B$ denotes the \emph{Cartesian product} of $A$ and $B$ and is defined as $\{ (a,b) \; ; \; a \in A,\ b \in B\}$. A \emph{binary relation}~\cite{bertrand1999new} $R$ defined on a set of arbitrary elements $X$ is a subset of $X \times X$, and we denote by $x \in R(y)$ or equivalently $x \  R \ y$ the fact that $(x,y) \in R$. An \emph{order relation}~\cite{bertrand1999new} is a binary relation $R$ which is reflexive, antisymmetric, and transitive. We denote by $\RC$ the binary relation on $X$ defined such that, $\forall x, y \in X$, $\left\{ x \ \RC \ y \right\} \Leftrightarrow \left\{ x \ R \ y \text{~and~} x \neq y \right\}$. A set $X$ of arbitrary elements supplied with an order relation $R$ on $X$ is called a \emph{partially ordered set} or \emph{poset} and is denoted by $(X,R)$, or shortly $|X|$ or even $X$ when no ambiguity is possible.

\medskip

\newcommand{\topology}{\mathcal{U}\xspace}

Let $X$ be a set of arbitrary elements, and let $\topology$ be a set of subsets of $X$. We say that $\topology$ is a \emph{topology} on $X$ if $\emptyset$ and $X$ are elements of $\topology$, if any union of elements of $\topology$ are elements of $\topology$, and if any finite intersection of elements of $\topology$ is an element of $\topology$. $X$ supplied with $\topology$ is denoted $(X,\topology)$ or shortly $X$ and is called a \emph{topological space}. The elements of $\topology$ are then called the \emph{open sets} of $X$ and any complement of an open set in $X$ is called a \emph{closed set} of $X$. We say that a subset of $X$ which contains an open set containing a point $x$ is a \emph{neighborhood} of $x$ in $X$. A topological space $X$ is said \emph{(topologically) connected} if it is not the disjoint union of two non-empty open sets.

\medskip

\newcommand{\RINVERSE}{R^{-1}\xspace}

As recalled in~\cite{boutry2020equivalence}, a \emph{$T_0$-space}~\cite{aleksandrov1956combinatorial,alexandroff2013topologie,kelley1955general}, let say $X$, is a topological space which satisfies the \emph{$T_0$ axiom of separation}: for two distinct elements $x,y$ of $X$, there exists a neighborhood of $x$ in $X$ which does not contain $y$ or a neighborhood of $y$ in $X$ which does not contain $x$. A \emph{discrete space}~\cite{alexandrov1937} is a topological space where any intersection of open sets is an open set. Posets are considered as topological spaces in the sense that we can induce a topology on any poset based on its order relation (Th. 6.52, p. 28 of~\cite{aleksandrov1956combinatorial}): for a poset $(X,R)$, the corresponding Alexandrov space is the topological space of domain $X$ where the closed sets are the sets $C \subseteq X$ such that $\forall x \in C$, $R(x)$ is included in $C$. Let us denote by $\RINVERSE$ the inverse of $R$. Then, by symmetry, we obtain that open sets are the sets $U$ such that for any $h \in U$, $\RINVERSE(h)$ is included into $U$. Discrete $T_0$-spaces are generally called \emph{Alexandrov spaces}. Details can be found in~\cite{eckhardt1994digital}.

\medskip

On an Alexandrov spaces $|X| = (X,R)$, for any element $h \in X$, we define respectively the \emph{(combinatorial) closure} of $h$:
$$\alpha(h) := \{h' \in X \; ; \; h' \in R(h)\},$$
its inverse operator called the \emph{(combinatorial) opening} of $h$:
$$\beta(h) := \{h' \in X \; ; \; h \in R(h')\},$$
and the \emph{neighborhood} of $h$:
$$\theta(h) := \{h' \in X \; ; \; h' \in R(h) \text{~or~} h \in R(h')\}.$$

Obviously, thanks to the properties explained before, for any $h \in X$, $\alpha(h)$ will be a closed set in $X$ and each $\beta(h)$ will be an open set in $X$. In other words, combinatorial and topological definitions are equivalent in Alexandrov spaces.

\medskip

The operators $\alpha$, $\beta$ and $\theta$ can also be defined for sets: $\forall S \subseteq X$, $\alpha(S) := \cup_{p \in S} \alpha_X(p)$, $\beta(S) := \cup_{p \in S} \beta_X(p)$, and  $\theta(S) := \cup_{p \in S} \theta_X(p)$, where $\alpha(S)$ is closed and $\beta(S)$ is open thanks to the properties exposed before.

\medskip

We call \emph{path}~\cite{bertrand1999new} in a poset $|X|$ a finite sequence $\langle p^0,\dots,p^k\rangle$ in $X$ when for all $i \in \INT{1}{k}$, $p^i \in \TC(p^{i-1})$. We say that a poset $|X|$ is \emph{path-connected}~\cite{bertrand1999new} if for any points $p,q$ in $X$, there exists a path into $X$ joining them. In posets, path-connectedness is equivalent to (topological) connectedness. The greatest path-connected set in the digital set $X$ containing $p \in X$ is called the \emph{connected component}~\cite{aleksandrov1956combinatorial} of $X$ containing $p$ and we denote it by $\COMPOS(X,p)$; by convention, when $p$ does not belong to $S$, we write $\COMPOS(S,p) = \emptyset$. Any non-empty subset of a poset $S$ which can be written $\COMPOS(S,p)$ for some $p \in S$ is called a \emph{connected component} of $S$. The set of connected components of a poset $S$ is denoted by $\COMPOS(S)$.

\medskip

Let $X$ be a poset of rank $n \geq 1$. We call \emph{$(n-1)$-path} of $\vert X \vert$ any path made of only $(n-1)$- and $n$-faces in $\vert X \vert$; two elements of $X$ joined by such a path are said to be $(n-1)$-connected in $\vert X \vert$. We say that this poset $\vert X \vert $ is \emph{$(n-1)$-connected} when any two $n$-faces of $\vert X \vert $ are $(n-1)$-connected.

\medskip

The \emph{rank} $\RANK(h,|X|)$ of an element $h$ in the poset $|X|$ is $0$ if $\AC_X(h) = \emptyset$ and is equal to:
$$\max_{x \in \AC_X(h)} (\RANK(x,|X|)) + 1$$
otherwise. The \emph{rank} of $|X|$ is denoted by $\RANK(|X|)$; it is equal to the maximal rank of its elements when it is not empty and $-1$ otherwise. An element $h$ of $X$ such that $\RANK(h,|X|) = k$ is called a \emph{$k$-face}~\cite{bertrand1999new} of $X$. 

\medskip

Let $|X| = (X, \alpha_X),|Y| = (Y, \alpha_Y)$ be two posets. We say that $|X|$ is a \emph{suborder} of $|Y|$ when $X \subseteq Y$ and $\alpha_X = \alpha_Y \cap X \times X$. In a natural way, we will have, for any $h \in X$, $\alpha_X(h) = \alpha_Y(h) \cap X$, $\beta_X(h) = \beta_Y(h) \cap X$, $\theta_X(h) = \theta_Y(h) \cap X$, $\AC_X(h) = \AC_Y(h) \cap X$, $\BC_X(h) = \BC_Y(h) \cap X$ and $\TC_X(h) = \TC_Y(h) \cap X$.

\medskip

In this paper, we will often say that a face $a \in X$ is \emph{parent} of a face $b \in X$ (relatively to the poset $\vert X \vert$) when $a \in \BC_X(b))$. In an equivalent manner and with the same hypothesis, we will say that $b$ is a \emph{child} of $a$. We will call a bijection $F : X \rightarrow Y$ from the poset $X$ to the poset $Y$ an \emph{(order) isomorphism} when $a \in \BC_X(b)$ is equivalent to $F(a) \in \BC_Y(F(b))$. Furthermore, we will say that two posets $X,Y$ are \emph{isomorphic} between them when their exist an isomorphism between them. 

\medskip

\begin{figure}
    \centering
    \includegraphics[width=0.5\linewidth]{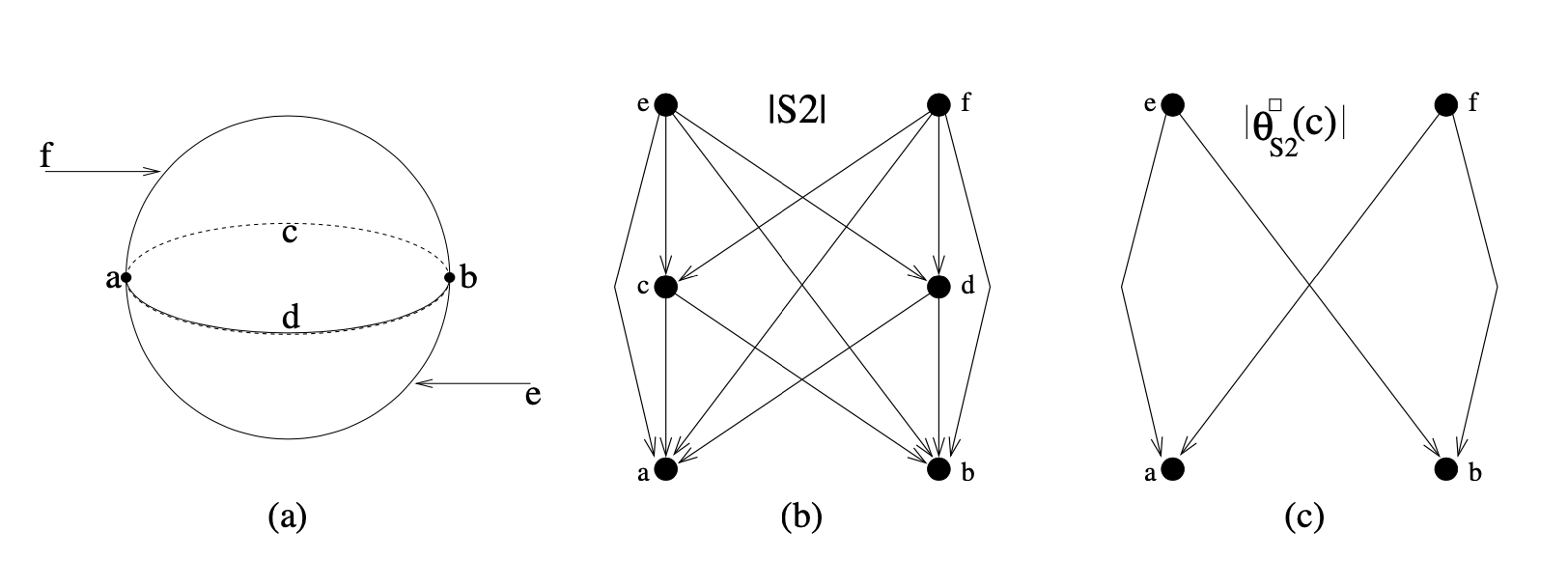}
    \caption{On the left side, an example of poset $2$-surface representing a sphere called $S2$. On its right side, its representation using a Hasse diagram with parenthood relationship encoded with arrows, and on the right side, the neighborhood $\TC_{S2}(c)$ of the face $c$ arbitrarily chosen in $S2$. Since this poset is connected and that for every face $h \in S2$, $\vert \TC_{S2}(c) \vert$ is a $1$-surface, then it is a discrete $2$-surface. This figure is extracted from~\cite{daragon2005surfaces}.}
    \label{fig:surface}
\end{figure}

Let $|X|$ be a suborder of $|Y|$. $|X|$ is said to be \emph{countable} if $X$ is countable. Also, $|X|$ is called \emph{locally finite} if for any element $h \in S$, the set $\theta_X(h)$ is finite. When $|X|$ is countable and locally finite, it is said to be a \emph{CF-order}~\cite{bertrand1999new} in $|Y|$.

\medskip

Let $|X|$ be a CF-order (in some poset $|Y|$); $|X|$ is said to be a (discrete) $(-1)$-surface if $X = \emptyset$, or a (discrete) $0$-surface if $X$ is made of two different faces $x,y \in X$ such that $x \not \in \TC_X(y)$, or a (discrete) $k$-surface, $k \in \llbracket 1,n\rrbracket$, if $|X|$ is connected and for any $h \in X$, $|\TC_X(h)|$ is a $(k-1)$-surface. According to Evako {\it et al.}~\cite{evako1996dimensional}, Khalimsky grids are discrete surfaces. Another example in a less constrained context is given in Figure~\ref{fig:surface}.

\medskip

A very strong property~\cite{daragon2005surfaces} of discrete surfaces is that the joint $X * Y$ is a $k+\ell+1$-surface \textbf{iff} $X$ is a $k$-surface and $Y$ an $\ell$-surface for some values $k, \ell \geq -1$.

\subsection{Recalls about borders and poset-based connected manifolds}

\newcommand{\border}{\Delta \xspace}
\newcommand{\Interior}{\mathrm{Int}\xspace}

\begin{figure}[h!]
    \centering
    \includegraphics[width=0.5\linewidth]{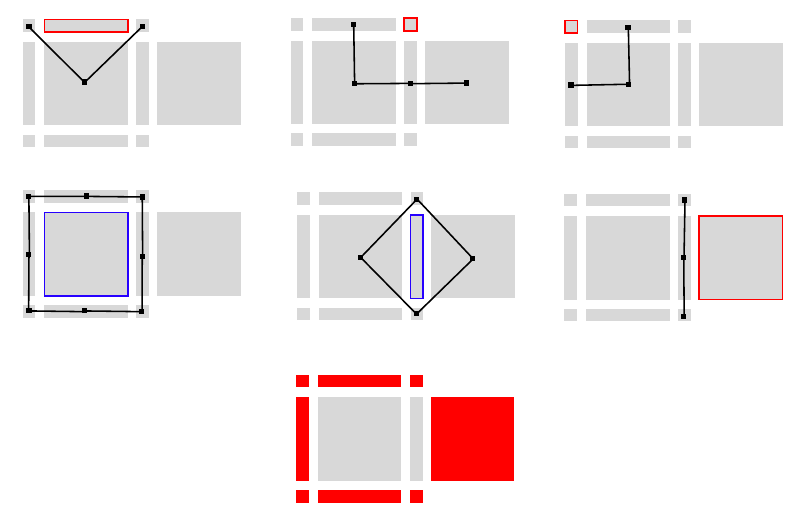}
    \caption{How we compute a border in a poset. In the raster scan order, we check whether the $1$-face encircled in red belongs to the border of the gray poset $X$. Since the neighborhood of this face is not a simple close curve (more exactly, it is not a discrete $1$-surface), then it belongs to $\border X$. We continue to follow the raster scan order to see that the reasoning is the same for the two $0$-faces encircled in red. At the second line, the two first faces (each one being encircled in blue) admit as neighborhoods a simple close curve, that is, a discrete $1$-surface, so these two faces belong to the interior of $X$. The last studied face on the second line has as neighborhood different from a $1$-surface, so this face belongs to the border $\border X$. The border is finally depicted in red at the third line. Note that our definition of border satisfies always that $\border X \subseteq X$, which is not always the case in topology (see the combinatorial boundaries as presented in~\cite{najman2013discrete}). This figure is extracted from~\cite{boutry2024introducing}.}
    \label{fig:border}
\end{figure}

\medskip

\begin{Definition}[Border and Interior~\cite{boutry2024introducing}]
Let $|X|$ be a poset of rank $n \geq 0$. We define the border of $|X|$ as the set:
$$\border X = \left\{ h \in X \; ; \; \left\vert \TC_X(h)\right\vert  \text{ is not a } (n - 1)\text{-surface}\right\}.$$
Let $X$ be a poset of rank $n \geq 0$. The \emph{interior} of $X$ is defined as the set:
$$\Interior(X) := X \setminus \border X.$$
\label{def.border}
\end{Definition}

An intuitive presentation of the computation of a border of a poset is given in Figure~\ref{fig:border}.

\medskip

\begin{Definition}[\cite{boutry2024introducing}]
We say that a poset $\vert X \vert$ of rank  $n \geq -1$ is \emph{coherent} when it is empty (case $n = -1$) or when for any $h \in X$, we have the following properties:
$$\left\{
\begin{array}{l}
\RANK\left(\left\vert \TC_X(h)\right\vert \right) = \RANK(\vert X\vert) - 1,\\
\text{and } \vert \TC_X(h)\vert  \text{ is coherent.}
\end{array}
\right.$$
\end{Definition}

For example, discrete $n$-surfaces are coherent posets.

\medskip

Let us now present a remarkable property of the border operator.

\begin{Proposition}[\cite{boutry2024introducing}]
Let $X$ be a coherent poset of rank $n \geq 0$, and $h$ be an element of $\border X$. Then we have the following remarkable property:
$$\border \TC_X(h) \subseteq \TC_{\border X}(h).$$
In other words, the border of the neighborhood is included in the neighborhood of the border.
\label{propo:border}
\end{Proposition}

\begin{figure}[h!]
    \centering
    \includegraphics[width=0.4\linewidth]{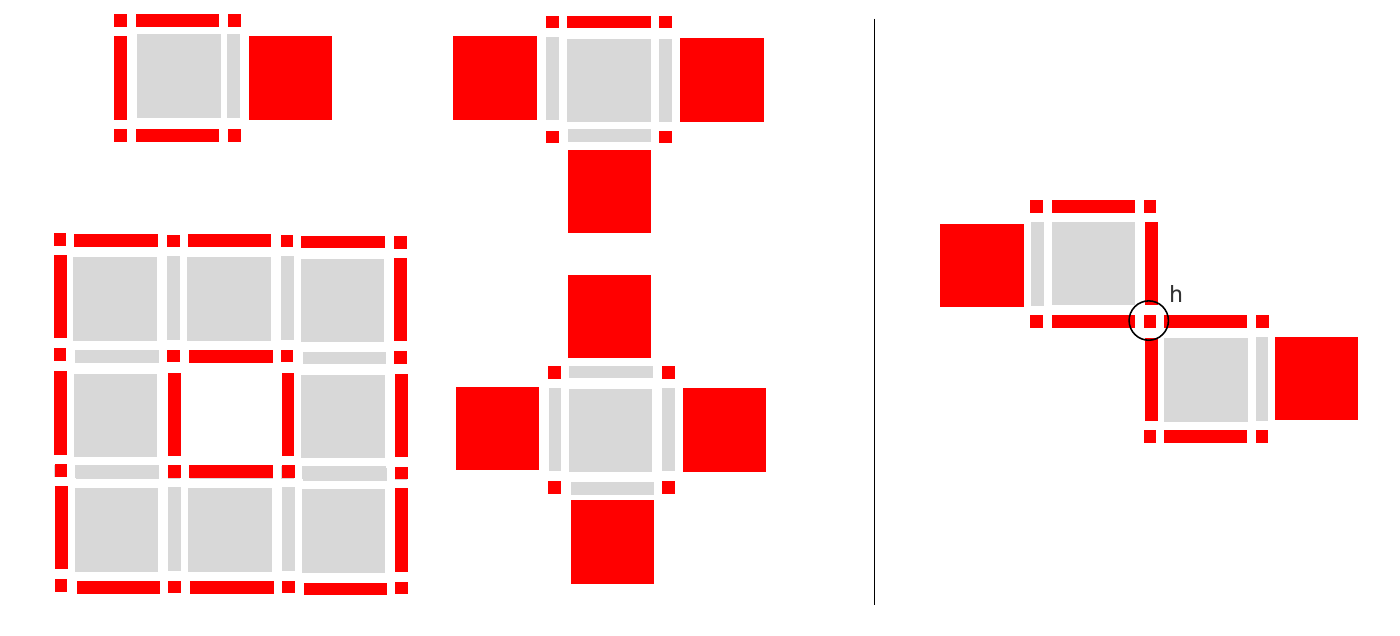}
    \caption{On the left side, several posets (the interior part is in gray and the border is in red) which are $2$-PCM's since the neighborhood of any of their face is connected and either a $1$-PCM or a $1$-surface. On the right side, a poset whose face $h$ satisfies that $\vert \TC_X(h) \vert$ is not connected and then not a $1$-surface, so this last poset in not a PCM (and $h$ is a pinch). This figure is extracted from~\cite{boutry2024introducing}.}
    \label{fig:pcms}
\end{figure}

\begin{Definition}[$n$-PCM~\cite{boutry2024introducing}]
Let $|X|$ be a poset of rank $n \geq - 1$. We say that $|X|$ is a $(-1)$-PCM when it is the empty order. We say that $|X|$ is a $0$-PCM when $|X|$ can be written $\{h\}$ with $h$ an arbitrary element; and we say that $|X|$ is an $n$-PCM, with $n \geq 1$, when $|X|$ is path-connected, $\border X \neq \emptyset$, and for any $h \in X$, either $\left\vert \TC_X(h)\right\vert $ is an $(n-1)$-surface (and $h$ does not belong to $\border X$), or it is an $(n-1)$-PCM (and $h$ belongs to $\border X$).
\end{Definition}

The difference between $n$-PCM's and \textquote{not $n$-PCM's} is depicted in Figure~\ref{fig:pcms}.

\medskip

We can strengthen the definition of an $n$-PCM by defining \emph{smooth} $n$-PCM's.

\begin{Definition}[Smooth $n$-PCM's~\cite{boutry2024introducing}]
Let $|X|$ be a poset of rank $n \geq -1$. We say that $|X|$ is a \emph{smooth $-1$-PCM} when $X = \emptyset$, that $|X|$ is a \emph{smooth $0$-PCM} when $X = \{h\}$ with $h$ some arbitrary element, and that $|X|$ is a \emph{smooth $n$-PCM}, $n \geq 1$, when it is connected, its border $\border X \neq \emptyset$ is a discrete $(n-1)$-surface or a separated union of discrete $(n-1)$-surfaces, and for any $h \in X$, $\vert \TC_X(h)\vert $ is either a smooth  $(n-1)$-PCM (border case) or an $(n-1)$-surface (interior case).
\end{Definition}

Naturally, $n$-PCM's and smooth $n$-PCM's are coherent.

\medskip

\begin{figure}[h!]
    \centering
    \includegraphics[width=0.7\linewidth]{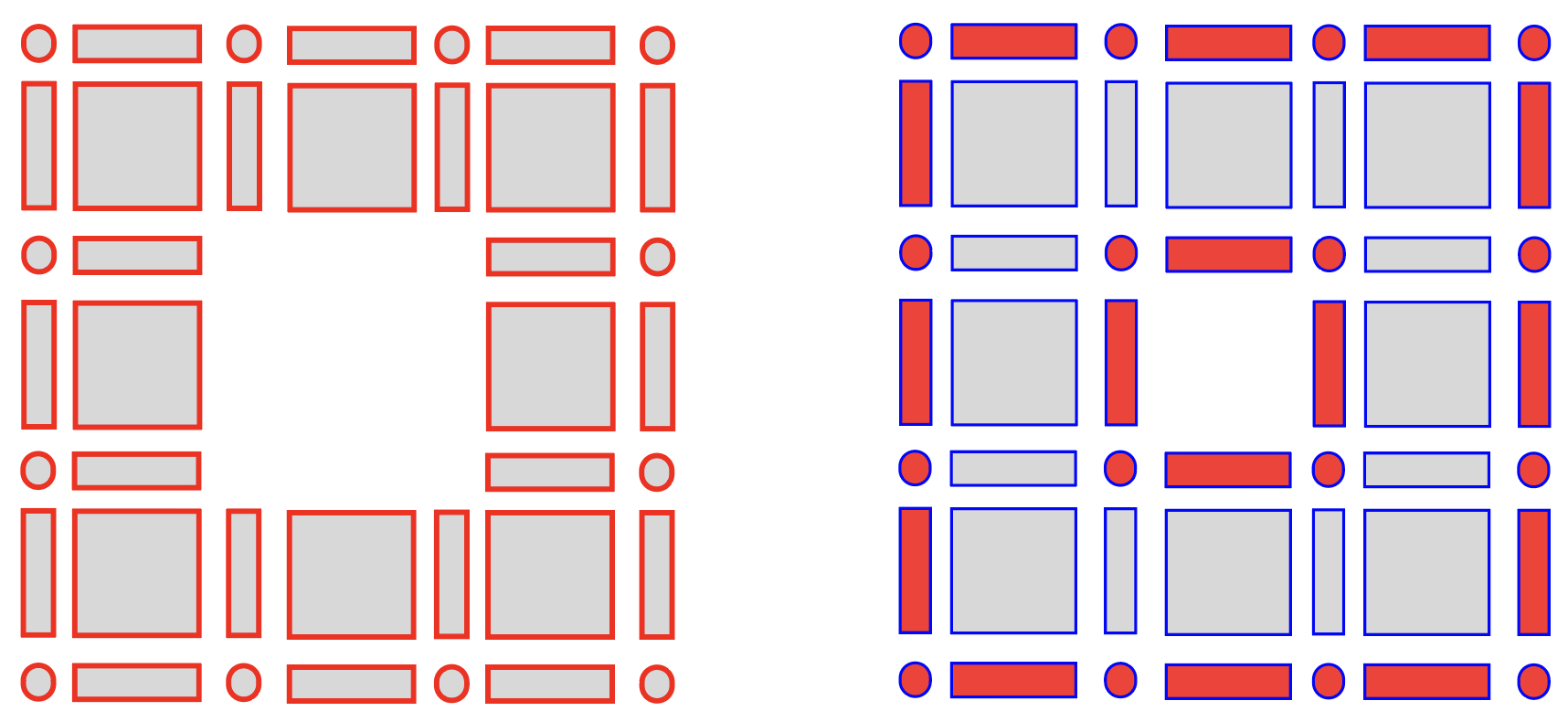}
    \caption{On the left side, an example of $2$-PCM $X$ which is not smooth, its border $\border X$ is made of the entire space ($\border X = X$) which is naturally not a $1$-surface. On the right side, a smooth $2$-PCM $Y$ whose border $\border Y$ is made of the separated union of two $1$-surfaces (depicted in red).}
    \label{fig:smooothORnot}
\end{figure}

We depict the difference between a $2$-PCM and a smooth $2$-PCM in Figure~\ref{fig:smooothORnot}.

\medskip

Now, let us recall briefly some important facts relative to (smooth) $n$-PCM's that have been shown in~\cite{boutry2024introducing}:

\newcommand{\family}{\mathcal{I}\xspace}
\newcommand{\ci}{C_i\xspace}

\begin{itemize}

    \item The only poset being at the same time a PCM and a discrete surface is the empty order. In other words, this set is isomorphic to a $-1$-surface and to a $-1$-PCM.

    \item Let $A,B$ be two (smooth) $n$-PCMs, with $n \geq 0$, such that $\RANK(A \cup B) = n$, $\border A = \border B = A \cap B$, and satisfying $\theta(\Interior(A)) \cap \Interior(B) =\emptyset$. Then, $A \cup B$ is a discrete $n$-surface. In other words, we can easily build a discrete $n$-surface starting from smooth $n$-PCM's when \textquote{natural} constraints are satisfied.

    \item Let $A,B$ be two (smooth) $n$-PCMs with $n \geq 1$. When $A \cap B = \border A \cap \border B$ is a (smooth) $(n-1)$-PCM, the rank of $A \cup B$ is equal to $n$, $\theta(A \setminus B) \cap (B \setminus A) = \emptyset$, and $\theta(\Interior(A)) \cap \Interior(B) = \emptyset$, then $A \cup B$ is an $n$-PCM. In other words, we can build $n$-PCM's starting from smooth $n$-PCM's using the union operator when some particular constraints are satisfied.

    \item Let $X$ be an $n$-surface with $n \geq 1$, and let $N$ be some suborder of $X$ which is an $(n-1)$-surface. We denote then by $\{\COMPOS_i\}_{i \in \family}$ the connected components of $\COMPOS(X \setminus N)$. We assume that the cardinality $\#(\family)$ of $\family$ is greater than or equal to $2$. Then, the components $\vert C_i \sqcup N\vert $ are (smooth) $n$-PCMs. In other words, we naturally obtain $n$-smooth PCM's by \textquote{cutting} a discrete $n$-surface.

    \item Let $n \geq 1$ be some integer. Let $X$ be a discrete $n$-surface and let $N \subset X$ be a discrete $(n-1)$-surface. Let us denote $\{\ci\}_{i \in \family} = \COMPOS(X \setminus N)$. Then, the cardinality of $\family$ is equal to or lower than two. In other words, the \textquote{cut} of a discrete surface using a discrete surface of lower rank leads to at most two (smooth) $n$-PCM's.
    
\end{itemize}

\subsection{Recalls about joins of surfaces and PCM's}

Let $\vert X\vert  := (X, R_X)$ and $\vert Y\vert  := (Y,R_Y)$ be two posets; it is said that $\vert X\vert $ and $\vert Y\vert $ can be \emph{joined}~\cite{bertrand1999new} if $X \cap Y = \emptyset$. If $\vert X\vert $ and $\vert Y\vert $ can be joined, the \emph{join} of $\vert X\vert $ and $\vert Y\vert $ is denoted $\vert X\vert * \vert Y\vert $ and is equal to: $$(X \cup Y, R_X \cup R_Y \cup X \times Y).$$

\medskip

Let $X,Y$ be two posets that can be joined. We recall the property from~\cite{daragon2005surfaces} that when $x \in X$, $\TC_{X * Y} = \TC_X(x) * Y$, and when $y \in Y$, $\TC_{X * Y}(y) = X * \TC_Y(y)$.

\medskip

\begin{Proposition}[\cite{boutry2024introducing}]
When $X$ is a $k$-surface and $Y$ an $\ell$-surface, then $X * Y$ is a $(k+\ell+1)$-surface. When $X$ is an $k$-PCM and $Y$ is an $\ell$-surface, then $X * Y$ and $Y * X$ are a $(k+\ell+1)$-PCM's. When $X$ is an $k$-PCM and $Y$ is an $\ell$-PCM, then then $X * Y$ is a $(k+\ell+1)$-PCM. 
\label{propo.joins}
\end{Proposition}

\begin{Proposition}[\cite{boutry2024introducing}]
The order joint of a $k$-PCM $P_k$ and of a $\ell$-surface $S_\ell$ is a $(k+\ell+1)$-PCM. Furthermore, $\border(P_k * S_\ell) = S_\ell \cup \border P_k$.
\label{propo.join.pcm.surf}
\end{Proposition}

Let us introduce a notation that will be very useful in the sequel.
\begin{Notation}
Let $|X|$ be a poset of rank $n \geq 0$. We denote the elements of rank $k \in \INT{0}{n}$ in $X$ by $(X)_k$. Furthermore, with $m,n$ two given integers s.t. $m \leq n$, we denote the set $$(X)_m \cup (X)_{m+1} \cup \dots \cup (X)_n$$
by $(X)_{[m,n]}$.
\end{Notation}

\subsection{Simplices and pseudomanifolds}

Let $\Lambda$ be a set of arbitrary elements. We call simplex every finite set $h \subseteq \Lambda$; its \emph{dimension}, denoted by $\DIM(h)$, is defined as its cardinality minus one. The simplices of cardinality $1$ (then of dimension $0$) are called \emph{vertices}.

\begin{Definition}
Let $X$ be a set of simplices. When for any $h \in X$,  any $h' \subseteq h$ belongs to $X$ too, we say that $X$ is closed under inclusion. When additionally for any $h,h' \in X$, the simplex $h \cap h'$ belongs to $X$ too, we say that $X$ is a simplicial complex. 
\end{Definition}

\begin{Definition}
Let $X$ be a simplicial complex of rank $n \geq 0$. When $n = 0$, $X$ is a degenerated set of one arbitrary element. When $n \geq 1$, then we say that $X$ is a \emph{pseudomanifold} when:
\begin{itemize}
    \item it is \emph{pure}, that is, for any $h \in X$, $\left(\beta_X(h)\right)_n \neq \emptyset$,
    \item each $(n-1)$-face of $X$ is the face of one \textbf{or} two $n$-faces of $X$,
    \item it is \emph{$(n-1)$-connected}, that is, any pair of $n$-faces of $X$ is joined by some path in $(X)_{[n-1,n]}$.
\end{itemize}
\end{Definition}

\begin{Definition}[\cite{daragon2005discrete}]
When two simplices $h,h'$ belong to the same simplicial complex $X$ with the properties $h \cap h' = \emptyset$, we will denote the \emph{simplicial join} $h \cup h'$ by $h \circ h'$.
\end{Definition}

\begin{Definition}
We recall that in a simplicial complex $X$ of rank $n \geq 0$, the \emph{link} of an element $h \in X$ is defined as:
$$\LINK(h,X)= \{ h' \in X \; ; \; h \circ h' \in X \}.$$
\end{Definition}

\begin{Definition}[\cite{daragon2005discrete}]
For two simplicial complexes $X,Y$ satisfying $X \cap Y = \emptyset$, we will define the \emph{simplicial join} $X \circ Y$ of $X$ and $Y$ as:
$$X \circ Y := X \cup Y \cup \left\{x \circ y \; ; \; x \in X, y \in Y\right\}.$$
Of course, a simplicial join of two simplicial complexes is a simplicial complex.
\end{Definition}

Among the different possible definitions~\cite{bertrand2023discrete,basak2020three,bagchi2023structure,padilla2002normal,novik2012face} of what is a normal pseudomanifold, we will chose this one:

\begin{Definition}
A pseudomanifold $X$ or rank $n \geq 0$ is called a \emph{\textbf{normal} pseudomanifold} when the link of every face of $X$ of dimension $d \leq n-2$ is a pseudomanifold.
\end{Definition}

The reason explaining we want that the links are pseudomanifolds is that we want they are connected~\cite{novik2012face}.

\medskip

We will also remark the difference with the definition proposed in~\cite{bertrand2023discrete} where any $(n-1)$-face has necessarily two parents, which implies that these structures do not have any border.

\begin{Proposition}[\cite{daragon2005surfaces}]
For a simplicial complex $X$ of rank $n \geq 0$, for any $h \in X$, $|\AC_X(h)|$ is a $(\DIM(h) - 1)$-surface. 
\label{propo.12.0}
\end{Proposition}

\medskip

We will call \emph{separated union} of two posets $X$ and $Y$ the the union $X \cup Y$ of the two disjoint posets, assuming that:
$$\theta_{X \cup Y}(X) \cap Y = \emptyset = \theta_{X \cup Y}(Y) \cap X,$$
which means that there exists no path joining any element of $X$ to any element of $Y$ in the poset $\vert X \cup Y \vert$.

\section{New properties relative to discrete surfaces, borders and PCM's in the context of simplicial complexes}
\label{sec.newbackground}
\begin{Proposition}
When $X$ is a simplicial complex of rank $n \geq 0$, then for any $x \in X$ and for any $y \in \AC_X(x)$, then $\left\vert \BC_X(y) \cap \AC_X(x)\right\vert$ is a $(\DIM(x) - \DIM(y) - 2)$-surface.
\label{propo.simplicial}
\end{Proposition}

\myproof We recall that when $X$ is a simplicial complex of rank $n \geq 0$, then $\vert \AC_X(x) \vert$ is a discrete $(\RANK(x,X)-1)$-surface by Proposition~\ref{propo.12.0}. Now, let $y$ be an element of $\AC_X(x)$, then $\left\vert \TC_{\AC_X(x)}(y)\right\vert$ is a discrete $(\RANK(x,X)-2)$-surface. However, 

$$\TC_{\AC_X(x)}(y) = \BC_{\AC_X(x)}(y) * \AC_{\AC_X(x)}(y) = (\BC_X(y) \cap \AC_X(x)) * \AC_X(y),$$

so $\vert \BC_X(y) \cap \AC_X(x) \vert$ is a $(\RANK(x,X) - \RANK(y,X) - 2)$-surface (since $\vert \AC_X(y) \vert $ is a discrete $(\RANK(y,X) - 1)$-surface). Since in simplicial complexes, the rank is equal to the dimension, we obtain that $\vert \BC_X(y) \cap \AC_X(x) \vert$ is a $(\DIM(x) - \DIM(y) - 2)$-surface.

\medskip

The proof is done. \qed

\medskip

Note that the previous property was firstly introduced for cubical grids~\cite{daragon2005surfaces,boutry2020equivalence}.

\subsection{Propositions relative to discrete surfaces and PCM's}

We recall that we know that discrete surfaces are pure posets~\cite{daragon2005discrete}. However the question relative to $n$-PCM's is not yet solved.

\begin{Proposition}
Any poset $X$ which is an $n$-PCM, $n \geq 0$ is pure.
\label{propo.PCMs.pure}
\end{Proposition}

\myproof Let us recall the definition of the concatenation operator. Let $X$ be some poset. For two paths $\pi_1 = \langle h_0,\dots,h_n\rangle$ of length $n \geq 0$ and $\pi_2 = \langle\ell_0,\dots,\ell_m\rangle$ of length $m \geq 0$, both in the poset $X$, and assuming that $h_n \in \TC_X(\ell_0)$, we can \emph{join} $\pi_1$ and $\pi_2$ using the \emph{concatenation operator} denoted by $\sim$ to make a new path in $X$:
$$\pi_1 \sim \pi_2 := \langle h_0, \dots,h_n,\ell_0,\dots,\ell_m\rangle$$
of length $n+m+1$ (the value $1$ is due to the concatenation operator).

\medskip
Let us now prove the assertion presented above. When $n=0$, every face is maximal in $X$ so this poset is pure. In the more general case $n \geq 1$, let us start from a $k$-face $h$ in $X$ with $k \in \INT{0}{n-1}$, otherwise $h$ is one more time already maximal. By definition of an $n$-PCM, we know that $\vert\TC_X(h)\vert = \vert\BC_X(h)\vert * \vert\AC_X(h)\vert$ is of rank $(n-1)$. Two cases are then possible: 

\begin{itemize}
    \item when $k=0$, $\AC_X(h) = \emptyset$, thus $\vert\BC_X(h)\vert$ is of rank $(n-1)$ and then we can find an increasing\footnote{We call a path \emph{increasing} relatively to $X$ when the rank of the faces of the path increases from the start to the end of the path.} path $\pi_\beta$ relatively to $X$ of length $(n-k-1) = (n-1)$ from a face $h_m$ of minimal rank in $\vert\BC_X(h)\vert$ to a face $h_M$ of maximal rank in $\vert\BC_X(h)\vert$. This way, by joining the degenerated path $\langle h \rangle$ to $\pi_\beta$, we obtain an increasing path (relatively to $X$) denoted by $\langle h \rangle \sim \pi_\beta$ which shows that $h_M$ is of rank $\RANK(\BC_X(h)) + 1 = n$. So, $h_M$ is maximal in $X$ and then $h$ admits a maximal face of rank $n$ in its neighborhood.

    \item when $k > 0$, we do the same procedure in $\vert\AC_X(h)\vert$ to obtain $\pi_\alpha$ from $h'_m$ (minimal in $\vert\AC_X(h)\vert$) to $h'_M$ (maximal in $\vert\AC_X(h)\vert$) of length $\RANK_X(h) - 1 = k - 1$, so the increasing final path:
    $$\Pi = \pi_\alpha \sim \langle h \rangle \sim \pi_\beta,$$
    which shows that $h_M$ is one more time maximal in $X$ since the length of the path $\Pi$ is equal to $(k-1) + (n - k - 1) + 2 = n$. Then, $h$ admits a maximal face of rank $n$ in its neighborhood.

\end{itemize}

The proof is done. \qed

\begin{Definition}
We say that a poset of rank $n$ is \emph{homogeneous} when for any element $h$ of $X$ of rank $k$ in $X$, there exists for any $\ell \in \INT{k+1}{n}$ at least one face of rank $\ell$ in $\vert\BC_X(h)\vert$ and for any $\ell \in \INT{0}{k-1}$ one face of rank $\ell$ in $\vert\AC_X(h)\vert$.
\end{Definition}

\begin{Proposition}
Let $X$ be an $n$-PCM, $n \geq -1$, then $X$ is homogeneous.
\label{propo.PCMs.homogeneous}
\end{Proposition}

\myproof The cases where $X$ is of rank $n \in \INT{-1}{0}$ are immediate. When $n \geq 1$, then for any element $h$ of rank $k \in \INT{0}{n}$, we know that there exists a path of $k$ elements of $X$ of increasing ranks $\langle h_0, \dots,h_\ell,\dots,h_{k-1} \rangle$ in $\vert\AC_X(h)\vert$ where the indices of each $h_\ell$ represent its rank in $X$. The first part of the assertion is proven.

\medskip

Now, we know from $\vert\TC_X(h)\vert = \vert\BC_X(h)\vert * \vert\AC_X(h)\vert$ that: 
$$\RANK(\vert\TC_X(h)\vert) = \RANK(\vert\BC_X(h)\vert) + \RANK(\vert\AC_X(h)\vert) + 1.$$
Since $\RANK(\vert\AC_X(h)\vert) = k - 1$ and $\RANK(\vert\TC_X(h)\vert) = n - 1$, then $\RANK(\vert\BC_X(h)\vert) = n - k - 1$. Thus there exists a path of $(n-k)$ elements of $X$ in $\BC_X(h)$ of increasing ranks $\langle h_{k+1}, \dots, h_m,\dots, h_n \rangle$. One more time, the indices of each $h_m$ represent its rank in $X$. The second part of the assertion is proven.

\medskip

The proof is done. \qed

\begin{Proposition}
When $X$ is an $n$-PCM, with $n \geq 0$, then $\vert \AC_X(h) \vert$ is a $(\RANK(h,X)-1)$-PCM or a $(\RANK(h,X)-1)$-surface.
\label{propo.PCMs.page1}
\end{Proposition}

\myproof When $n = 0$, then $h \in X$ is of rank $0$ in $X$ and satisfies that $\AC_X(h)$ is a $-1$-surface (and at the same time a $-1$-PCM). When $n \geq 1$, let us proceed by decreasing ranks of $h$ in $X$:
\begin{itemize}
    
    \item When $h_n$ is of rank $n$ in $X$, $\vert\AC_X(h_n)\vert = \vert\TC_X(h_n)\vert$ which is a $(n-1)$-PCM or a discrete $(n-1)$-surface by definition of an $n$-PCM.

    \item When $h_{n-1}$ is of rank $(n-1)$ in $X$, it admits a parent $h_n$ in $X$ by purity of $n$-PCM (see Proposition~\ref{propo.PCMs.pure}). Since $h_{n-1} \in \AC_X(h_n)$, we can write:
    
    $$\TC_{\AC_{X}(h_n)}(h_{n-1}) = \BC_{\AC_{X}(h_n)}(h_{n-1}) * \AC_{\AC_{X}(h_n)}(h_{n-1}) = \AC_{\AC_{X}(h_n)}(h_{n-1}) = \AC_X(h_{n-1}),$$
    so $\AC_X(h_{n-1})$ is equal to $\TC_{\AC_{X}(h_n)}(h_{n-1})$ which is an $(n-2)$-PCM or a discrete $(n-2)$-surface since $\vert\AC_X(h_n)\vert$ is a $(n-1)$-PCM or a discrete $(n-1)$-surface.

    \item By homogeneity of an $n$-PCM (see Proposition~\ref{propo.PCMs.homogeneous}), we can continue this way for each $h_k$ of rank from $(n-2)$ to $0$ by choosing a parent $h_{k+1}$ of rank $(k+1)$ in $X$. The property $\BC_{\AC_{X}(h_{k+1})}(h_{k}) = \emptyset$ will lead to the desired property.
    
\end{itemize}

The proof is done. \qed

\begin{Proposition}
When $X$ is an $n$-PCM, with $n \geq 0$, then $\vert \BC_X(h) \vert$ is a $(n - \RANK(h,X)-1)$-PCM or a $(n - \RANK(h,X)-1)$-surface.
\label{propo.PCMs.page2}
\end{Proposition}

\myproof The proof follows the same procedure as the one described in the previous property, except that we start from faces of rank $0$ to finish with the faces of rank $n$. The proof is done. \qed

\begin{Proposition}
Let $X$ be some $n$-PCM, with $n \geq 0$. Then, for any $h \in X$,
$$\vert\TC_X(h)\vert = \vert\BC_X(h)\vert * \vert\AC_X(h)\vert$$
with:
\begin{itemize}
    \item When $\vert\TC_X(h)\vert$ is a discrete $(n-1)$-surface,
        \begin{itemize}
        \item $\vert\BC_X(h)\vert$ is a discrete $(n-\RANK(h,X) -1)$-surface,
        \item $\vert\AC_X(h)\vert$ is a discrete $(\RANK(h,X) -1)$-surface,
        \end{itemize}
    \item When $\vert\TC_X(h)\vert$ is an $(n-1)$-PCM,
        \begin{itemize}
        \item $\vert\BC_X(h)\vert$ is a discrete $(n-\RANK(h,X) -1)$-surface or an $(n-\RANK(h,X) -1)$-PCM,
        \item $\vert\AC_X(h)\vert$ is a discrete $(\RANK(h,X) -1)$-surface or an $(\RANK(h,X) -1)$-PCM.
        \end{itemize}
\end{itemize}
\label{propo.PCMs.page3}
\end{Proposition}

\myproof This proof is simply the merge of Propositions~\ref{propo.PCMs.page1} and~\ref{propo.PCMs.page2}, associated with the fact that the joint of two discrete surfaces is a discrete surface. The proof is done. \qed

\begin{Proposition}
Let $X$ be an $n$-PCM, with $n \geq 1$, and let $a,b$ be two elements of $X$ satisfying $a \in \BC_X(b)$. Then,  $\vert \AC_X(a) \cap \BC_X(b)\vert$ is of rank $r = (\RANK(a,X) - \RANK(b,X) - 2)$ and it is either a discrete $r$-surface of a $r$-PCM.
\label{propo.alpha/inter.star.in.PCMs}
\end{Proposition}

\myproof The term $\vert \AC_X(a) \cap \BC_X(b)\vert$ can be rewritten $\vert \BC_{\AC_X(a)}(b)\vert$ with $\AC_X(a)$ either a $(\RANK(a,X) - 1)$-surface or a $(\RANK(a,X) - 1)$-PCM (by Proposition~\ref{propo.PCMs.page3}). So, still by Proposition~\ref{propo.PCMs.page3}, $\vert \AC_X(a) \cap \BC_X(b)\vert$ is of rank equal to
$$\RANK(\vert \AC_X(a) \vert) - \RANK(b,\vert \AC_X(a) \vert) - 1 = \RANK(a,X) - 1 - \RANK(b, X) - 1 = \RANK(a,X) - \RANK(b, X) - 2,$$
and is either a PCM or discrete surface. The proof is done. \qed

\subsection{Properties of the border (Part 1)}

\begin{Lemma}
Let $X$ be a simplicial complex of rank $n \geq 0$, and $h$ be an element of $\border X$. Then we have the following remarkable property:
$$\TC_{\border X}(h) \subseteq \border \TC_X(h).$$
\label{Lemma:border.inverse}
\end{Lemma}

\myproof Let $X$ be some simplicial complex of rank $n \geq 0$. Let $h$ be some element of $X$, and let $h'$ be some element of $\TC_X(h)$.

\medskip

We obtain naturally that $h' \in \border X$ and that $h' \in \TC_X(h)$.

\medskip

We want to show that $h' \in \border \TC_X(h)$, that is, that:

$$\left\vert \TC_{\TC_X(h)}(h') \right\vert $$

is \emph{not} a discrete $(n-2)$-surface.

\medskip

\underline{First case :} when $h' \in \BC_X(h)$, we can decompose this term in this way:

$$\left\vert \TC_{\TC_X(h)}(h') \right\vert 
 = \left\vert \BC_{\TC_X(h)}(h') \right\vert * \left\vert \AC_{\TC_X(h)}(h') \right\vert  = \left\vert \BC_X(h') \right\vert * \left\vert \AC_{\TC_X(h)}(h') \right\vert, $$

where we can decompose the last term:

$$\left\vert \AC_{\TC_X(h)}(h') \right\vert = \vert \AC_X(h') \cap \TC_X(h) \vert = \left\vert \TC_{\AC_X(h')}(h) \right\vert = \left \vert \BC_{\AC_X(h')}(h) \right\vert * \left\vert \AC_{\AC_X(h')}(h) \right\vert,$$
equal to $\left\vert \BC_{X}(h) \cap \AC_X(h')\right\vert  * \left\vert \AC_X(h)\right\vert$
which is the join of the $(\RANK(h',X) - \RANK(h,X) - 2)$ by Proposition~\ref{propo.simplicial} and the $(\RANK(h,X)-1)$-surface, that is, a $(\RANK(h',X) - 2)$-surface. Now, let us assume that $\left\vert \BC_X(h') \right\vert$ (H) is a $(n - \RANK(h',X) - 1)$-surface, it would imply that $$\vert \TC_X(h') \vert = \vert \BC_X(h') \vert * \vert \AC_X(h') \vert$$
is a discrete $(n - \RANK(h',X) - 1) + (\RANK(h',X) - 1) + 1 = (n-1)$-surface, which is impossible since $h' \in \border X$. So (H) is false. This implies that $\left\vert \BC_X(h') \right\vert$ is not a discrete $(n - \RANK(h',X) - 1)$-surface, so $h' \in \border \TC_X(h)$.

\medskip

\underline{Secund case :} when $h' \in \AC_X(h)$, we obtain a similar reasoning.

\medskip

Thus, we have in any case that $\TC_{\border X}(h) \subseteq \border \TC_X(h)$ and the proof is done. \qed

\medskip

This leads to the strong following proposition.

\begin{Proposition}
Let $X$ be a \emph{coherent} simplicial complex of rank $n \geq 0$, and $h$ be an element of $\border X$. Then we have the following remarkable property:
$$\TC_{\border X}(h) = \border \TC_X(h).$$
\label{propo:border.neighborhood}
\end{Proposition}

\myproof This assertion is the consequence of Lemma~\ref{Lemma:border.inverse} and Proposition~\ref{propo:border}.\qed

\medskip

\begin{Proposition}
Let $X$ be some simplicial complex. When $X$ is an $n$-PCM with $n \geq 1$, then, for any $h \in \border X$, we have the following properties:
\begin{itemize}
    \item $\AC_X(h) \subseteq \border X$, that is, $\border X$ is closed under inclusion, and
    \item $\border \BC_X(h) \subseteq \border X$, that is, the border of the order made of the parents of $h$ in $X$, is included in the border of $X$.
\end{itemize}
\label{propo.bord.key}
\end{Proposition}

\myproof Let $h$ be an element of $\border X$. Let us study the term $\border \TC_X(h)$ in two ways. First,
$$\border \TC_X(h) = \border (\BC_X(h) * \AC_X(h)) = \border (\BC_X(h)) \cup \AC_X(h)$$
by Proposition~\ref{propo.join.pcm.surf}. Secund, using Proposition~\ref{propo:border.neighborhood}:
$$\border \TC_X(h) = \TC_{\border X}(h) = \TC_X(h) \cap \border X \subseteq \border X.$$
The union of these two results leads to:
$$\border (\BC_X(h)) \cup \AC_X(h) \subseteq \border X,$$
which concludes the proof. \qed

\subsection{Properties relative to PCM's that are simplicial complexes}

\begin{Proposition}
Let $X$ be at the same time a simplicial complex of rank $n \geq 1$ and an $n$-PCM. Then we have the following remarkable property for any $h \in \border X$:
$$\BC_{\border X}(h) = \border \BC_X(h).$$
\label{propo:border.neighborhood.beta}
\end{Proposition}

\myproof Let us choose some integer $n \geq 1$. When $h$ is a $0$-face, we can deduce by Proposition~\ref{propo:border.neighborhood} that:
$\BC_{\border X}(h) = \TC_{\border X}(h) = \border \TC_X(h) = \border \BC_X(h)$. Now, let $h$ be an element of $\border X$ of rank in $\INT{1}{n-1}$ and let us study in two different ways the term $\vert \TC_{\border X}(h)\vert$. First,
$$\vert \TC_{\border X}(h)\vert = \vert \BC_{\border X}(h) \vert * \vert \AC_{\border X}(h)\vert = \vert \BC_{\border X}(h) \vert * \vert \AC_X(h)\vert.$$
Secund, we have thanks to Proposition~\ref{propo:border.neighborhood} and Proposition~\ref{propo.join.pcm.surf}:
$$\vert \TC_{\border X}(h)\vert = \vert \border \TC_X(h) \vert = \vert\border(\BC_X(h) * \AC_X(h))\vert = \vert\border(\BC_X(h))\vert * \vert\AC_X(h)\vert.$$
The two expressions are joints with the same non-empty term $\vert \AC_X(h) \vert$, so $\vert\border \BC_X(h)\vert = \vert\BC_{\border X}(h)\vert$. The proof is done. \qed

\begin{Proposition}
Let $X$ be at the same time a simplicial complex of rank $n \geq 1$ and an $n$-PCM. Then we have the following remarkable property:
$$\LINK(h,\border X) \simeq \border \LINK(h,K).$$
\label{propo:border.neighborhood.link}
\end{Proposition}

\myproof Let us choose some integer $n \geq 1$. We assume that $h$ is an element of the poset $\border X$ where $X$ is at the same time a simplicial complex and a $n$-PCM. Since $\border X$ is a simplicial complex, $\LINK(h,\border X) \simeq \BC_{\border X}(h)$. Using Proposition~\ref{propo:border.neighborhood.beta}, we know that $\BC_{\border X}(h) = \border \BC_X(h)$, and since $X$ is a simplicial complex (by hypothesis), $\BC_X(h) \simeq \LINK(h,X)$ which implies that $\border \BC_X(h) \simeq \border \LINK(h,X)$. By grouping together all these isomorphisms, the proof is done. \qed

\begin{Proposition}
When $X$ is at the same time a simplicial complex of rank $n \geq 0$ and an $n$-PCM, then its border $\border X$ is of rank $(n-1)$.
\label{propo.border.natural}
\end{Proposition}

\myproof Let us treat the case $n = 0$. For any $h \in X$, $\vert \TC_X(h) \vert$ is the empty order, that is, a $-1$-surface. So $\border X = \emptyset$ and then has a rank $-1$. 

\medskip

Now, let us treat the case $n \geq 1$. We are going to proceed by induction on the rank $n$ of $X$ to prove that the rank of $\border X$ is at least equal to $(n-1)$. We define $(H_n)$ as \textquote{For any simplicial complex $X$ which is an $n$-PCM with $n \geq 1$, $\border X$ contains at least a face of rank $(n-1)$}. 

\medskip

\underline{Initialization $(n = 1)$:} when $X$ is at the same time a simplicial complex and a $1$-PCM, it is a simple open path $\langle h^0, h^1,\dots,h^L\rangle$, $L \geq 0$, with $h^0$ and $h^L$ two $0$-faces. Since $\vert \TC_X(h^0) \vert = \vert \{h^1\}\vert $ is a $0$-PCM, and then not a $0$-surface, $h^0$ belongs to $\border X$. Thus, $(H_1)$ is true.

\medskip

Now we assume that the properties $(H_1)$ until $(H_{n-1})$ are true, and we want prove that it implies that $(H_n)$ is true.

\medskip

\underline{Heredity $(n \geq 2)$:} Since $X$ is an $n$-PCM with $n \geq 1$, then $\border X \neq \emptyset$. Since, by Proposition~\ref{propo.bord.key}, $\border X$ is closed under inclusion, then it contains at least one face $h_0$ of rank $0$ in $X$ (otherwise $\border X$ would be empty). Since $h_0 \in \border X$ with $X$ an $n$-PCM, then $\vert \TC_X(h_0) \vert$ is a $(n-1)$-PCM. Since $\TC_X(h_0) = \BC_X(h_0)$ which is isomorphic to the simplicial complex $\LINK(h_0,X)$, $\LINK(h_0,X)$ satisfies the conditions of $(H_{n-1})$ and contains at least one face $h_{n-2}$ of rank $(n-2)$ in $X$. Since $h_{n-2}$ is in the link of $h_0$ in $X$, the face $h_0 \circ h_{n-2}$ belongs to $\BC_X(h_0)$. However, $h_0 \circ h_{n-2}$ is the simplicial join of a simplex of dimension $(n-2)$ with a vertex, it is then a simplex of dimension $(n-1)$. Consequently, there exists at least one $(n-1)$-face in $\BC_X(h_0)$. Since $\BC_X(h_0) \subseteq \border X$ by Proposition~\ref{propo.bord.key}, $(H_n)$ is then true.

\medskip

By induction on $n$, we have proven that $\border X$ is at least of rank $(n-1)$.

\medskip

Now, let us prove (by contradiction) that, for $X$ a simplicial complex and $n$-PCM of rank $n \geq 1$, the rank of $\border X$ is at most of rank $(n-1)$. Let us assume that $\border X$ contains some $n$-face $h_n$ of $X$. The direct consequence is that $\vert \TC_X(h_n) \vert = \vert \BC_X(h_n) \vert * \vert \AC_X(h_n) \vert = \vert \AC_X(h_n)\vert $ is a $(n-1)$-PCM since $X$ is an $n$-PCM. However, $X$ is also a simplicial complex, so $\vert \TC_X(h_n)\vert  = \vert \AC_X(h_n)\vert $ is a $(n-1)$-surface. Then, $\vert \TC_X(h_n) \vert$ it at the same time a $(n-1)$-surface and an $(n-1)$-PCM with $(n-1)) \geq 0$, we have a contradiction. So the rank of $\border X$ is at most $(n-1)$.

\medskip

By grouping the two proven properties, the rank of $\border X$ is equal to $(n-1)$. The proof is done. \qed

\medskip

\begin{Lemma}
When $|X|$ is an $n$-PCM, $n \geq 0$, and when for some $h \in X$, $|\TC_X(h)|$ is the joint of a poset $A$ with another poset $B$, then when $A$ is a $k$-surface, $B$ is a $(n - k - 1)$-PCM, and when $B$ is a $\ell$-surface, $A$ is a $(n - \ell -1)$-PCM.
\label{lemme.star.2024}
\end{Lemma}

\myproof We can proceed by induction on the hypothesis $(H_n)$ defined as: \textquote{when $X,Y$ are two posets such that their join $\vert X\vert  * \vert Y\vert $ is an $n$-PCM, $n \geq 0$, with $\vert X\vert $ an $k$-surface, $k \geq 0$, then $\vert Y\vert $ is an $(n - k - 1)$-PCM.}.

\medskip

\underline{Initialization $(n = 1)$:} when $\vert X\vert  * \vert Y\vert $ is a $0$-PCM, then $\vert X\vert $ is necessarily a $0$-surface and $\vert Y\vert $ the empty order, that is, a $(-1)$-PCM. $(H_1)$ is then true.

\medskip

\underline{Heredity ($n \geq 2$):} we assume that $(H_m)$ is true from $m = 0$ to $m = (n-1)$, let us prove it for the case $n$. When $\vert X\vert  * \vert Y\vert $ is an $n$-PCM, we obtain that for any $x \in X$, $|\TC_{X * Y}(x)| = |\TC_X(x) * Y|$ is an $(n-1)$-surface (Case $A$) or an $(n-1)$-PCM (Case $B$) by definition of an $n$-PCM. Let us assume that we are in Case $A$. Since $\vert X\vert $ is a $k$-surface, then $|\TC_X(x)|$ is a $(k-1)$-surface, from which we deduce that $\vert Y\vert $ is a $((n-1) - 1 - (k-1) = (n-k-1)$-surface (since $|\TC_X(x) * Y|$ is an $(n-1)$-surface). It would imply that $\vert X\vert  * \vert Y\vert $ is a $n$-surface, which is impossible by hypothesis (the only poset being at the same time a discrete surface and a PCM is the empty order and here $n \geq 1$). Then we are in Case $B$ and $|\TC_{X * Y}(x)|$ is an $(n-1)$-PCM. Using the induction hypothesis $(H_{n-1})$ on the join $\vert \TC_X(x)\vert  * \vert Y\vert $ which is a $(n-1)$-PCM with $|\TC_X(x)|$ an $(k-1)$-surface, we obtain that $\vert Y\vert $ is a $((n-1) - (k-1) - 1) = (n-k-1)$-PCM.

\medskip

The second part of the proposition follows a symmetrical reasoning. The proof is done. \qed

\begin{Proposition}
Let $\vert X\vert $ be at the same time an $n$-PCM with $n \geq 0$ and a simplicial complex, then, when $h \in \triangle X$, $\vert \BC_X(h)\vert $ is a $(n - \DIM(h) - 1)$-PCM, and  when $h \not \in \triangle X$, $\vert \BC_X(h)\vert $ is a $(n - \DIM(h) - 1)$-surface.
\label{propo.12}
\end{Proposition}

\myproof When $h$ does not belong to the border $\border X$, $|\TC_X(h)|$ is an $(n-1)$-surface, and since $|\TC_X(h)| = |\BC_X(x)| * |\AC_X(x)|$, the two terms ate discrete surfaces. Since $|\AC_X(x)|$ has its rank equal to the rank of $x$ in $X$, and then equal to the dimension of $x$ (since we are in a simplicial complex), the rank of $|\BC_X(x)|$ is equal to $(n - \DIM(h) - 1)$.

\medskip

Now, when $h$ belongs to the $\border X$, then $|\TC_X(h)|$ is an $(n-1)$-PCM (by definition of an $n$-PCM). Since $|\TC_X(h)| = |\BC_X(x)| * |\AC_X(x)|$ with $|\AC_X(h)|$ a $(\DIM(h) - 1)$-surface (by Proposition~\ref{propo.12.0}), we obtain that $|\BC_X(x)|$ is of rank $(n-\DIM(h)-1)$ and by Lemma~\ref{lemme.star.2024} that it is a PCM. The proof is done. \qed

\subsection{Properties of the border (Part 2)}
\begin{Proposition}
When $X$ is at the same time a simplicial complex of rank $n \geq 0$ and an $n$-PCM, then its border $\border X$ is made of the closure (by inclusion) of its $(n-1)$-simplices:
$$\border X = \left\vert\left\{h \; ; \; h \in \alpha_X(h_{n-1}) \; , \; h_{n-1} \in (\border X)_{n-1}\right\}\right\vert$$
\label{propo.border.closure.NM1.simplices}
\end{Proposition}

\myproof Since we already know that $X$ is a simplicial complex, we just want to prove that for any face $h \in \border X$, there exists some $(n-1)$-face $h_{n-1} \in \border X$ satisfying the relation $h \in \alpha_X(h_{n-1})$.

\medskip

We recall that no $n$-simplex can belong to $\border X$. For any $(n-1)$-face $h \in \border X$, $h \in \alpha_X(h)$ (by reflexivity of a order relation), so the property is satisfied. Then we can only treat the case where $h \in \border X$ with $\DIM(h) \leq n - 2$.

\medskip

Since $h$ is an element of $\border X$, then $\vert \TC_X(h)\vert$ is not a discrete $(n-1)$-surface. Because $X$ is a $n$-PCM, $\vert\TC_X(h)\vert$ is an $(n-1)$-PCM. Since $X$ is a simplicial complex, $\vert\AC_X(h)\vert$ is a $(\DIM(h)-1)$-surface, thus $\vert \BC_X(h) \vert$ is a $(n - \DIM(h) - 1)$-PCM by Lemma~\ref{lemme.star.2024}. By Proposition~\ref{propo.border.natural}, the border of $\vert\LINK(h,X)\vert$ (which is at the same time a simplicial complex and an $(n - \DIM(h) - 1)$-PCM) is of rank $(n - \DIM(h) - 2)$, thus $\vert\border \BC_X(h)\vert$ is of rank $(n - \DIM(h) - 2)$ too. In other words, by Proposition~\ref{propo:border.neighborhood.beta}, $\vert\BC_{\border X}(h)\vert$ is of rank $(n - \DIM(h) - 2)$, which means that there exists at east one parent $h'$ of $h$ in $\BC_{\border X}(h)$ of dimension:
$$\DIM(h) + 1 + (n - \DIM(h) - 2) = (n-1).$$
Since $\BC_{\border X}(h) = \border \BC_X(h) \subseteq \border X$ by Proposition~\ref{propo.bord.key}, $h'$ belongs to $\border X$. The proof is done. \qed

\section{Relating pseudomanifolds, discrete surfaces, and PCM's}
\label{sec.equivalence}

Let us now present some strong relations that exist between these topological structures.

\subsection{From discrete $n$-surfaces and $n$-PCM's to normal pseudomanifolds}

In this subsection, we are going to study which relation to normal pseudomanifolds we obtain when we assume that we have a simplicial complex which is either a discrete $n$-surface or an $n$-PCM.

\begin{figure}[h!]
    \centering
    \includegraphics[width=0.35\linewidth]{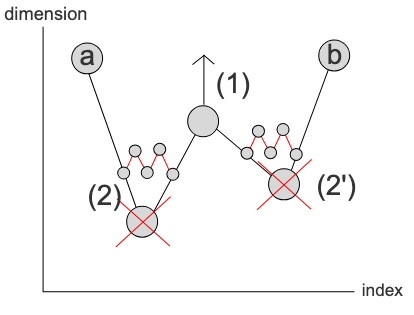}
    \caption{The procedure we propose to find a $(n-1)$-path between the two $n$-faces $a,b$ of $X$, made only (at the end of the procedure) of $(n-1)$- and $n$-faces of $X$. This procedure is possible whenever $X$ is a $n$-PCM or a discrete $n$-surface and shows that $X$ is $(n-1)$-connected. Step 1: starting from any path going from $a$ to $b$, first we deduce a new equivalent path (in the sense same start and same end) whose even faces along the path are $n$-faces. Step 2: starting from this new path, and thanks to the properties of discrete $n$-surfaces and of $n$-PCM's to satisfy that $\vert\BC_X(h)\vert$ is a non-empty connected subset of $X$ as soon as the rank of $h$ is X is strictly lower than $(n-1)$, we can replace each face $p^{2k+1}$ with a minimal rank with the path in its open neighborhood $\vert\BC_X(p^{2k+1})\vert$ of greater minimal rank. At the end, we have the path we looked forward, which proves that $X$ is $(n-1)$-connected.}
    \label{fig:procedure}
\end{figure}

\begin{Proposition}
Let $X$ be a simplicial complex of rank $n \geq 1$. When $\vert X\vert $ is a $n$-PCM or an $n$-surface, it is $(n-1)$-connected. Furthermore, for any $h \in X$ with $\DIM(h) \in \INT{1}{n-2}$, $\vert \BC_X(h)\vert $ is $(n-1)$-connected too.
\label{propo.nmoinsunconnexe}
\end{Proposition}

\myproof The reasoning of the proof of this proposition is detailed in Figure~\ref{fig:procedure}.

\medskip

Let $a_n, b_n$ be two $n$-faces of $X$. By definition of a $n$-PCM, $X$ is connected since $n \geq 1$. Then, there exists a path in $X$ joining $a_n$ and $b_n$.

\medskip

We can assume without any constraint that this path is of the form: $\pi = \langle p^0 = a_n, \dots, p^{2 K} = b_n \rangle$ with for any $k \in \INT{0}{K-1}, p^{2 k} \in \BC_X(p^{2 k + 1})$ and for any $k \in \INT{1}{K}$, $p^{2 k - 1} \in \AC_X(p^{2 k})$.

\medskip

Now, since $X$ is pure, for any $k \in \INT{1}{K-1}$, for any term of the form $p^{2 k}$, we can find a new face $q_n^{2 k} \in \left(\beta_X(p^{2 k})\right)_n$. Since we have:
$$\left\{p^{2 k - 1},p^{2 k + 1}\right\} \subseteq \AC_X(p^{2 k}),$$
this face will satisfy:
$$\left\{p^{2 k - 1},p^{2 k + 1}\right\} \subseteq \AC_X(q_n^{2 k}),$$
from which we can deduce a new path $\pi' = \langle a_n, p^1, q_n^2, \dots, p^{2 K -3}, q_n^{2K-2}, p^{2 K - 1}, b_n \rangle$ joining $a_n$ and $b_n$ in $X$.

\medskip

Now, if any odd term in $\pi'$ belongs already to $(X)_{n-1}$, $\pi'$ is an $(n-1)$-path joining $a_n$ and $b_n$ in $X$ and the proof is done. However, if there exists any term $p^{2 k + 1}$ in $\pi'$ which does not belong to $(X)_{n-1}$, it belongs to $(X)_{[0,n-2]}$ and we can replace it by a new sequence where the minimal rank will be greater to the one of $p^{2 k + 1}$.

\medskip

Indeed, let us localize the element $p^{2 k + 1}$ in $\pi'$ of minimal rank in $\pi'$, and let us denote $r = \RANK(p^{2k+1},X) \in [0,n-2]$ its rank in $X$. Then, $\vert \BC_X(p^{2k+1})\vert $ is, by Propositions~\ref{propo.12.0} and~\ref{propo.12}, either a $(n - r - 1)$-surface or a $(n-r-1)$-PCM: since $r \leq n - 2$, then $(n - r - 1) \geq 1$ and thus $\vert \BC_X(p^{2k+1})\vert $ is connected. The two terms $q_n^{2k}$ and $q_n^{2k+2}$ around $p^{2k+1}$ in the path $\pi'$ belonging to $\BC_X(p^{2k+1})$, there exists then a path $\PILOCAL$ in $(X)_{[r+1,n]}$ joining $q_n^{2k}$ and $q_n^{2k+2}$. We can then replace the sequence $\langle q_n^{2k},p^{2k+1},q_n^{2k+2}\rangle$ with $\PILOCAL$ in $\pi'$.

\medskip

Since $X$ is finite, we can repeat this procedure until there is not anymore elements of rank in $[0,n-2]$ in the studied path, and this way we obtain a final path $\pi^*$ joining $a_n$ and $b_n$ in $(X)_{[n-1,n]}$, that is, a $(n-1)$-path. So, $a_n$ and $b_n$ are $(n-1)$-connected. We have proven that $\vert X \vert$ is $(n-1)$-connected.

\medskip

Let us now treat the second part of the assertion. By following a similar procedure as before, we can easily prove that any two $n$-faces $a_n$ and $b_n$ belonging to $\BC_X(h)$ will be connected by a $(n-1)$-path $\pi$ subset of $\BC_X(h)$. Indeed, let $p^{2k+1}$ be one of the faces of $\pi$ whose dimension is minimal (among the faces of $\pi$). If its dimension is $(n-1)$, $\pi$ is an $(n-1)$-path which connects $a_n$ and $b_n$ and there is nothing more to prove. However, when its dimension is $r \leq n - 2$, $\vert\BC_X(p^{2k+1})\vert$ is connected and thus we can replace in $\pi$ the sequence $\langle p^{2k}, p^{2k+1}, p^{2k+2} \rangle$ with a new path $$\langle p^{2k}, \dots, p^{2k+2} \rangle \subseteq \BC_X(p^{2k+1}) \subseteq \BC_X(h)$$ where the face(s) of minimal degree will have a dimension greater than $r$. By proceeding in this way for all the faces in $\pi$ of dimension $r$, we will obtain a new path which is at the same time in $\BC_X(h)$ and in $(X)_{[r+1,n]}$. We repeat this procedure until $\pi$ belongs to $(X)_{[n-1,n]}$, always by preserving $\pi \in \BC_X(h)$, which is ensured to be possible since $X$ is finite. The proof is done. \qed 

\begin{Lemma}
Let $X$ be a simplicial complex of rank $n \geq 1$. When $\vert X\vert $ is a $n$-PCM or an $n$-surface, it is a pseudomanifold of rank $n$.
\label{lemma.pseudomanifold.sens.simple}
\end{Lemma}

\myproof Since $\vert X\vert $ is either a discrete $n$-surface of an $n$-PCM, it is pure. Then, to prove that $\vert X\vert $ is a pseudomanifold of rank $n$, we have to prove $(A)$ that each $(n-1)$-face of $\vert X\vert $ is the face of one or two faces of $\vert X\vert $, and $(B)$ that $\vert X\vert $ is $(n-1)$-connected.

\medskip

\textbf{(A):} Let $h_{n-1}$ be some $(n-1)$-face of $\vert X\vert $. When $\vert X \vert$ is an $n$-PCM, when $h \in \triangle X$, then $\vert \BC_X(h_{n-1})\vert $ is by Proposition~\ref{propo.12} a $(n - \DIM(h_{n-1}) - 1) = 0$-PCM, that is, it is a (degenerated) set of one $n$-face of $X$. When $h \not \in \triangle X$, then $\vert\BC_X(h_{n-1})\vert$ is a $0$-surface by Proposition~\ref{propo.12.0} ($n$-surface case) and by Proposition~\ref{propo.12} ($n$-PCM case), thus $h$ is a face of two $n$-faces of $X$.

\medskip

\textbf{(B):} Since $\vert X\vert$ is either a discrete $n$-surface or an $n$-PCM with $n \geq 1$, it is $(n-1)$-connected by Proposition~\ref{propo.nmoinsunconnexe}.

\medskip

The proof is done. \qed

\begin{Theorem}
Let $X$ be a pure simplicial complex of rank $n \geq 2$. When $\vert X\vert $ is an $n$-PCM or a discrete $n$-surface, it is a normal pseudomanifold of rank $n$.
\label{theorem.un}
\end{Theorem}

\myproof By Lemma~\ref{lemma.pseudomanifold.sens.simple}, $X$ is a normal pseudomanifold when for any $h \in X$ with $\DIM(h) \leq n - 2$, $\vert \LINK(h,X)\vert $ is a pseudomanifold (of rank $(n - \DIM(h) - 1)$).

\medskip

Let us denote by $d$ the dimension of $h$ and by $r$ the rank of $\vert \LINK(h,X)\vert $. We want to prove that $(A)$ $\vert \LINK(h,X)\vert $ is of rank $r = (n - d - 1)$, that $(B)$ $\vert \LINK(h,X)\vert $ is a simplicial complex (see the two conditions hereafter), $(C)$ $\vert \LINK(h,X)\vert $ is pure, $(D)$ that each $(r-1)$-face of $\vert \LINK(h,X)\vert $ is the face of one or two $r$-faces of $\LINK(h,X)$, and $(E)$ that for each pair of facets $a,b$ in $\LINK(h,X)$, there exists a path in $\vert (\LINK(h,X))_{[r-1,r]}\vert $ joining $a$ and $b$.

\medskip

\begin{figure}
    \centering
    \includegraphics[width=0.25\linewidth]{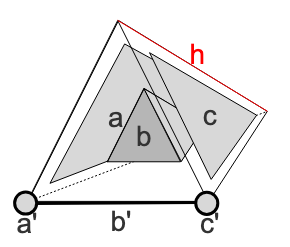}
    \caption{Presentation of the isomorphism between $\BC_X(h)$ and $\LINK(h,X)$ in a simple case: we easily observe that $b' = b \circ h$ is the common parent of $a' = a \circ h$ and $b' = b \circ h$ if and only if $b$ is the common parent of $a$ and $c$.}
    \label{fig:isomorphism}
\end{figure}

Let $h$ be some element of $X$ with $\DIM(h) \leq n - 2$. We define the following mapping: $$\MAPPING : \LINK(h,X) \rightarrow X : h' \rightarrow \MAPPING(h') = h \circ h'.$$
As seen before, this mapping is an isomorphism (see Figure~\ref{fig:isomorphism}) from $\vert\LINK(h,X)\vert$ to $\vert\BC_X(h)\vert$.

\medskip

\textbf{(A):} Since $\vert \LINK(h,X)\vert $ and $\vert \BC_X(h)\vert $ are isomorphic, then $\vert \LINK(h,X)\vert $ is of rank $r = n - d - 1$ by Propositions~\ref{propo.12.0} and~\ref{propo.12} whatever $h$ belongs to $\triangle X$ or not.

\medskip

\textbf{(B):} The first condition to prove is that $\LINK(h,X)$ is closed under inclusion. Let $a$ be a face of $\LINK(h,X)$, then $a \circ h \in X$. Now, let $b \in X$ be a face of $a$. Since $a \circ h$ belongs to $X$ which is closed under inclusion, $b \circ h \in X$, thus $b$ belongs to $\LINK(h,X)$. Now, let us prove that when $a,b$ are two faces of $\LINK(h,X)$, their intersection belongs to $\LINK(h,X)$ too. Let $a,b$ be two faces of $\LINK(h,X)$ where, without constraint, $a \cap b \neq \emptyset$. Then, the faces $a \circ h$ and $b \circ h$ belong to $X$. Since $(a \cap b) \circ h = (a \circ h) \cap (b \circ h)$ which belongs to the simplicial complex $X$, then $(a \cap b) \circ h \in X$, that is, $a \cap b$ belongs to $\LINK(h,X)$.

\medskip

\textbf{(C):} We want to prove that $\vert \LINK(h,X)\vert $ is pure, that is, for any $h' \in \LINK(h,X)$, there exists a face of rank $r$ relatively to the suborder $\vert \LINK(h,X)(h')\vert $ and in $\vert \beta_{\LINK(h,X)}(h')\vert $. Let $h'$ be some face of $\LINK(h,K)$. Then, starting from $h'$, we can compute $\MAPPING(h') = h \circ h' \in \BC_X(h)$. Since $X$ is pure by hypothesis, there exists an $n$-face $h_n \in \BC_X(h \cup h')$ of maximal dimension $n$ in $\BC_X(h)$. From $h_n$, we can deduce by the inverse isomorphism $\MAPPING^{-1}$ a face $h_r = \MAPPING^{-1}(h_n) = h_n \setminus h$ of maximal dimension $r$ in $\vert \LINK(h,X)\vert$. Since $h_n \in \beta_X( h\circ h')$, thus $h_r \in \beta_{\LINK_X(h,X)}(h')$. Thus, $\vert \LINK(h,X)\vert $ is pure.

\medskip

\textbf{(D):} We want to prove that each $(r-1)$-face of $\vert \LINK(h,X)\vert $ is the face of one or two $r$-faces of $\vert \LINK(h,X)\vert$. We do almost the same procedure as before using the isomorphism $\MAPPING$: let $h'$ be some face of rank $(r-1)$ in $\LINK(h,X)$. Using $\MAPPING$, we obtain $h \circ h'$ into $\BC_X(h)$ and of dimension $(n-1)$. Then, $h \circ h'$ is the face of either one or two $n$-faces due to the fact that $X$ is a pseudomanifold. Using $\MAPPING^{-1}$, we obtain the faces in $\vert \LINK(h,X)\vert $ of rank $r$ the face $h'$ is the face of, whose cardinality is equal to $1$ or $2$.

\medskip

\textbf{(E):} Now, we want to prove that for each pair of maximal faces $a,b$ in $\vert \LINK(h,X)\vert $, there exists a path in $\vert (\LINK(h,X)_{[r-1,r]})\vert $ joining $a$ and $b$. For this aim, we simply start from the two faces $a_r,b_r$ of maximal rank $r$ in $\vert \LINK(h,X)\vert $, we compute their mappings by $\MAPPING$, $a_n$ and $b_n$ respectively, in $\BC_X(h)$. Since $\DIM(h) \leq n - 2$ by hypothesis, $\vert \BC_X(h)\vert $ is $(n-1)$-connected by Proposition~\ref{propo.nmoinsunconnexe} so there exists an $(n-1)$-path $\pi$ joining $a_n$ and $b_n$ in $\vert \BC_X(h)\vert $. Using now $\MAPPING^{-1}$ on $\pi$, we directly obtain a $(r-1)$-path in $\vert \LINK(h,X) \vert$ joining $a_r$ and $b_r$.

\medskip

The proof is done. \qed

\subsection{From normal pseudomanifolds to $n$-surfaces and $n$-PCM's}

Let us announce whether the converse implication is true.

\begin{Theorem}
Let $X$ be a normal pseudomanifold of rank $n \geq 0$. Then, either $\triangle X = \emptyset$ and $X$ is an $n$-surface, or $\triangle X \neq \emptyset$ and $X$ is an $n$-PCM.
\label{th.deux}
\end{Theorem}

\myproof Let $X$ be a normal pseudomanifold of rank $n \geq 0$. When $n = 0$, the pseudomanifold $X$ is a degenerated set made of one vertex, then it is a $0$-PCM. When $n \geq 1$ and $\triangle X = \emptyset$, then $\vert X\vert $ is connected (since it is a pseudomanifold) and additionally, for any $h \in X$, $\vert \TC_X(h)\vert $ is an $(n-1)$-surface, so $X$ is a $1$-surface.

\medskip

The case that we have to treat now is then when $X$ admits a (non-empty) border ($\triangle X \neq \emptyset$) with $n \geq 1$. Let $\vert X\vert $ satisfy these properties and let us prove it is an $n$-PCM. $\vert X\vert $ is a pseudomanifold, then it is connected. Also, for any $h \in X$, when $h \not \in \triangle X$, $\vert \TC_X(h)\vert $ is an $(n-1)$-surface. We want then to show that, when $h \in \triangle X$, $\vert \TC_X(h)\vert $ is an $(n-1)$-PCM.

\medskip

Let us prove that by induction with as induction hypothesis for $k \geq 1$:
$$H_k := \left\{ \forall \ h \in \triangle X, \ \vert \BC_X(h)\vert  \text{ is an $(n - \DIM(h)-1)$-PCM with $\DIM(h) = n - k$}\right\}.$$

\textbf{Initialization $(k = 0)$:} When $\DIM(h) = n$, we obtain naturally that for any face $h$ of dimension $n$, $\vert \BC_X(h) \vert$ is the empty order, that is, a $(-1)$-PCM and $H_0$ is true.

\medskip

\textbf{Initialization $(k = 1)$:} When $\DIM(h) = n - 1$, then $\BC_X(h)$ is made of one or two $n$-faces since $X$ is a pseudomanifold. However, the case with two elements, that is, two $n$-faces of $X$, is impossible since in this case we obtain that $\vert \BC_X(h)\vert $ is a $0$-surface, and then by the joining operation, $\vert \TC_X(h)\vert $ is an $(n-1)$-surface, thus $h \not \in \triangle X$, which contradicts the hypotheses. Thus, $\vert \BC_X(h)\vert $ is a $0$-PCM (made of only one $n$-face) and $H_1$ is true.

\medskip

\textbf{Heredity ($k \geq 2$)}: we assume that $H_{k'}$ is true for $k' \in [0,k-1]$, and we want to prove that this implies $H_k$. Now, let us treat the case $\DIM(h) = n - k \geq 2$. Since $\vert \BC_X(h)\vert $ is isomorphic to $\vert \LINK(h,X)\vert $ and $X$ is normal by hypothesis, then $\vert \BC_X(h)\vert $ is connected. Now, we want to prove that for any $h' \in \BC_X(h)$, the term $\left\vert \TC_{\BC_X(h)}(h')\right\vert $ is either an $(n-\DIM(h)-2)$-surface or an $(n-\DIM(h)-2)$-PCM.

\begin{itemize}
    \item When $h' \not \in \triangle X$, then $\vert \BC_X(h')\vert $ is an $(n - \DIM(h') - 1)$-surface, so we obtain:
    $$\vert \TC_{\BC_X(h)}(h')\vert  = \vert \BC_X(h')\vert  * \vert (\AC_X(h') \cap \BC_X(h))\vert $$
    which is the join of a $(n - \DIM(h') - 1)$-surface by Proposition~\ref{propo.12.0} and a $(\DIM(h') - \DIM(h) - 2)$-surface by Proposition~\ref{propo.simplicial}, it is then a $(n - \DIM(h) - 2)$-surface by Proposition~\ref{propo.joins}.

    \item When $h' \in \triangle X$: since $h' \in \BC_X(h)$, we have $\DIM(h') \geq \DIM(h)+1$, that is, $$\DIM(h') \in \INT{n - k + 1}{n} = \INT{n - (k - 1)}{n - 0},$$ which corresponds to the hypotheses $H_{k-1}$, ..., $H_{0}$, assumed to be true (they are induction hypotheses). Thus, we know that $\vert \BC_X(h')\vert $ is an $(n - \DIM(h') - 1)$-PCM, so we obtain:
    $$\left\vert \TC_{\BC_X(h)}(h') \right\vert  = \vert \BC_X(h') \vert 
 * \vert (\AC_X(h') \cap \BC_X(h))\vert $$
    which is the join of a $(n - \DIM(h') - 1)$-PCM by Proposition~\ref{propo.12.0} and a $(\DIM(h') - \DIM(h) - 2)$-surface by Proposition~\ref{propo.simplicial}, it is then a $(n - \DIM(h) - 2)$-PCM.

\end{itemize}

The proof is done. \qed

\subsection{Generalization of our main results}

Let us conclude this section with an important theorem.

\begin{Theorem}
When $X$ is a simplicial complex of rank $n \geq 0$,
\begin{itemize}
    \item either $\triangle X = \emptyset$, and in this case $X$ is a discrete surface iff it is a normal pseudomanifold,
    \item or $\triangle X \neq \emptyset$, and $X$ is a $n$-PCM  iff it is a normal pseudomanifold.
\end{itemize}
In other words, $n$-PCM's and discrete surfaces are a \textbf{generalization} of normal pseudomanifolds to posets.
\end{Theorem}

\section{PCM's vs. smooth PCM's}
\label{sec.pcm.vs.smooth.pcms}

Let us study whether $n$-PCM's which are supplied with additional topological properties (like being a simplicial complex) implies that it is a smooth $n$-PCM.

\begin{Theorem}
An $n$-PCM , $n \geq 0$, which is a simplicial complex is \textbf{not} necessarily a smooth $n$-PCM.
\end{Theorem}

\begin{figure}[h!]
    \centering
    \includegraphics[width=0.45\linewidth]{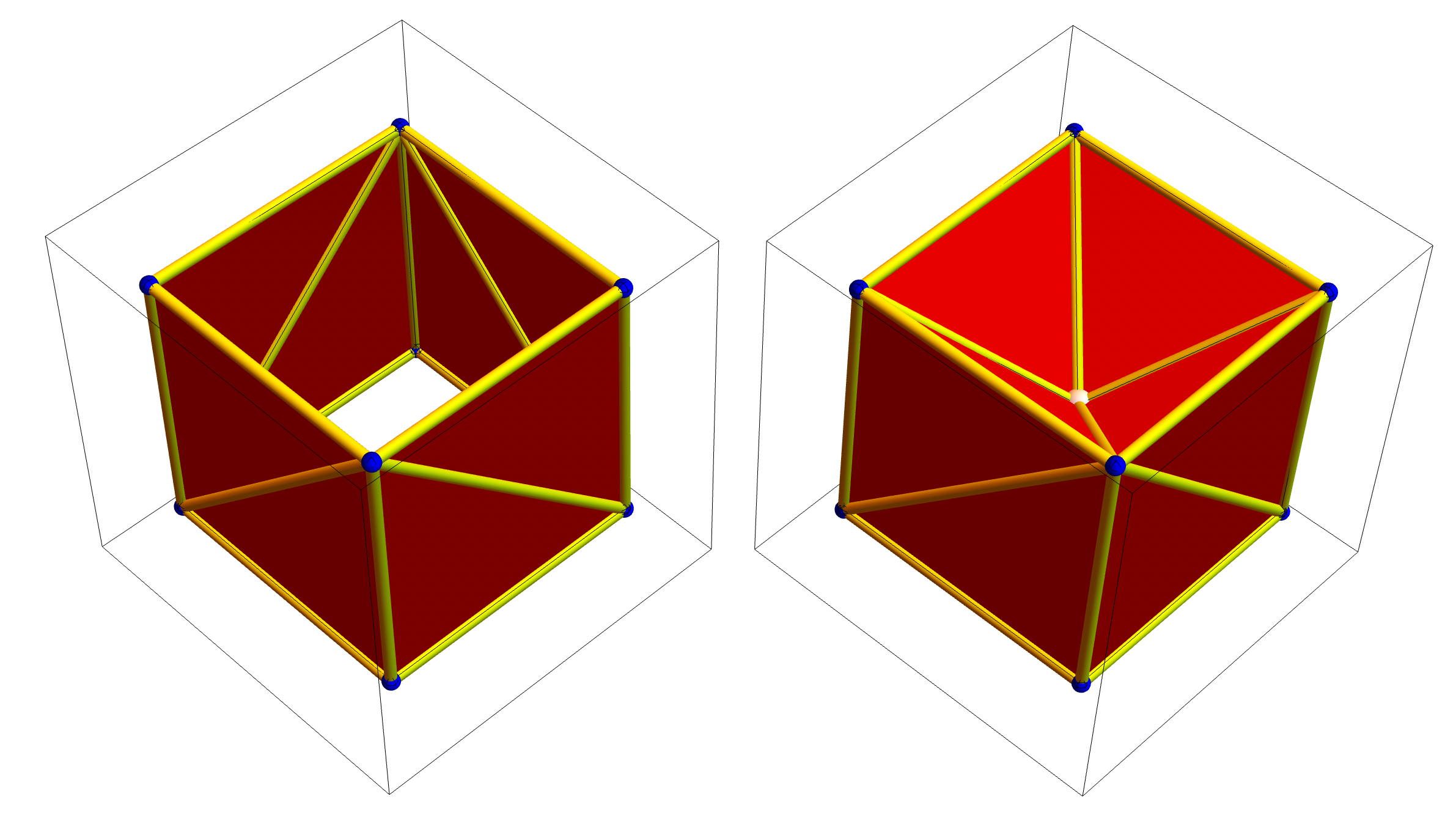}
    \caption{A 3D counter-example: on the right side, a simplicial complex $X$ equal to the simplicial join of two simplicial complexes: a $2$-PCM (on the left side and depicting a simplicial band) and a $0$-PCM $\vert \{h\}\vert$ (the vertex $h$ colored in white is placed at the barycenter of $Y$); the simplicial complex $X$ \textbf{is} a $3$-PCM but \textbf{not} a smooth $3$-PCM (as detailed in the body of the paper). We call this combinatorial structure $X$ the \emph{pinched simplicial box}.}
    \label{fig:contre-ex}
\end{figure}
\begin{figure}[h!]
    \centering
    \includegraphics[width=0.45\linewidth]{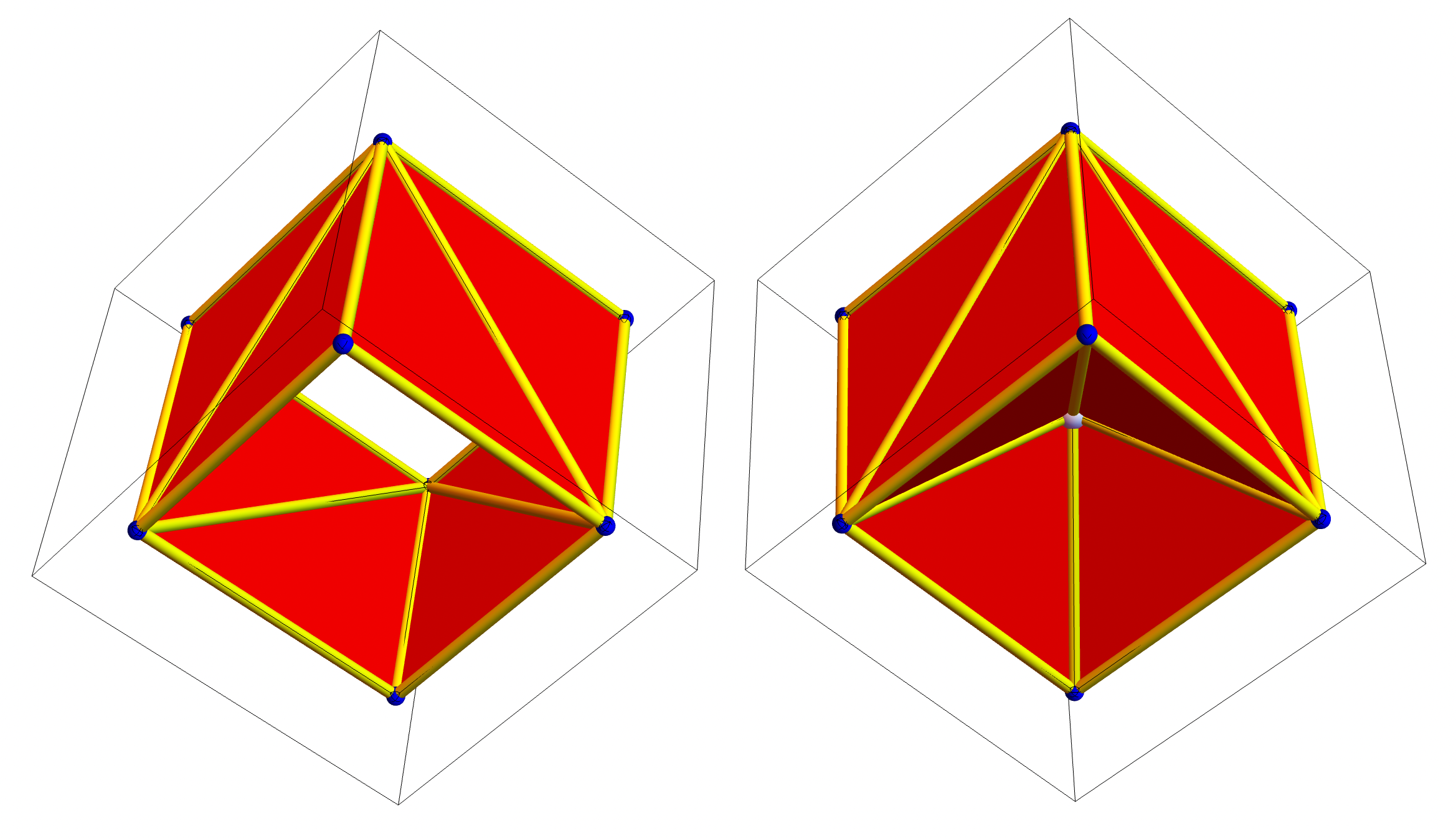}
    \caption{The same counter-example as in the previous figure, but seen from downside.}
    \label{fig:contre-ex-p2}
\end{figure}

\myproof Let us construct a poset $X$ by computing the simplicial join of a $2$-PCM $Y$ (on the left side in Figures~\ref{fig:contre-ex} and~\ref{fig:contre-ex-p2}) and a $0$-PCM $\vert\{h\}\vert$ made of the vertex $h$ (placed at its barycenter). We obtain a simplicial complex $X$ which is a $3$-PCM since it is connected, its border $\border X$ is non-empty, and at every element $h' \in X$, $\vert\TC_X(h')\vert$ is a discrete $2$-surface or a $2$-PCM (even when $h' = h$ since its link $\vert \LINK(h,X) \vert$ in $X$, isomorphic to $\vert \TC_X(h) \vert$, is the $2$-PCM $Y$). However, it is not a smooth $3$-PCM since $\vert \LINK(h,\border X)\vert \simeq \vert \BC_{\border X}(h)\vert = \vert \TC_{\border X}(h)\vert$ is made of two separated $1$-surfaces and thus it is not a discrete $1$-surface (the separated union of two discrete $1$-surfaces is not a discrete $1$-surface since it is not connected). The proof is done. \qed

\medskip

Let us assert that in the context of (coherent) simplicial complexes and that assuming a particular constraint, $n$-PCM's are also smooth $n$-PCM's.

\begin{Theorem}
Let $X$ be some $n$-PCM.
\begin{itemize}
    
    \item When $n \in \INT{-1}{0}$, then it is also a smooth $n$-PCM.

    \item When $X$ is a simplicial complex of rank $n = 1$, then $X$ is also a smooth $n$-PCM.
    
    \item When $X$ is a simplicial complex of rank $n \geq 2$, and when the constraint
    $$\forall \ h \in \border X, \vert \TC_{\border X}(h)\vert \text{ is a discrete $(n - 2)$-surface}, \ \ \ \ (C)$$
    is satisfied, then $X$ is also a smooth $n$-PCM. 
    
\end{itemize}
\label{theorem:smooth.under.condition}
\end{Theorem}

\myproof Let $X$ be an $n$-PCM of rank $n \geq -1$. We want to prove that $X$ is also a smooth $n$-PCM.

\medskip

\underline{Case $n = -1$:} $X$ is the empty order when $n = - 1$, which is also a smooth $-1$-PCM.

\medskip

\underline{Case $n = 0$:} $X$ has the form $X = \vert \{h\} \vert$ with $h$ some arbitrary element, and it is then also a smooth $0$-PCM.

\medskip

Now, we are going to proceed by induction on the rank $n \geq 1$ of $X$: let $(H_n)$ be the property that \textquote{Any simplicial complex $X$ which is a $n$-PCM and which satisfies $(C)$ when $n \geq 2$ is a smooth $n$-PCM}.

\medskip

\underline{Initialization (case $n = 1$):} Since $X$ is a simplicial complex and at a same time a $1$-PCM, it is a simple open path of the form $\langle h_0,\dots,h_{2K} \rangle $ with $K \geq 0$ an integer, $\{h_{2 k}\}_{k \in \INT{0}{K}}$ a set of vertices, and $\{h_{2 k - 1}\}_{k \in \INT{1}{K}}$ a set of $1$-simplices. Let us prove this is a smooth $1$-PCM. First, it is connected. Secund, its border is not empty ($\border X = \vert \{h_0,h_{2K}\}\vert$). Third, for any $h \in X$, $\vert \TC_X(h) \vert$ is either a discrete $0$-surface or a smooth $0$-PCM. Fourth, the border of $X$ is a $0$-surface. $(H_1)$ is then true.

\medskip

\underline{Heredity (case $n \geq 2$):} Now, assuming that $(H_k)$ is true for $k \in \INT{1}{n-1}$, let us prove that $(H_n)$ is true. Let $X$ be some $n$-PCM. First, since $X$ is a $n$-PCM with $n \geq 0$, it is connected. Secund, $X$ is an $n$-PCM with $n \geq 1$, then $\border X$ is non-empty. Third, we want to show that for any $h \in X$, $\vert\TC_X(h)\vert$ is either a smooth $(n-1)$-PCM or a discrete $(n-1)$-surface. When $h \in \Interior(X)$, $\vert\TC_X(h)\vert$ is a discrete $(n-1)$-surface by definition of the interior. When $h \in \border X$, $\vert\TC_X(h)\vert$ is a $(n-1)$-PCM. To prove that it is a smooth $(n-1)$-PCM, we can use the constraint $(C)$ which asserts that $\vert\TC_{\border X}(h)\vert$ is a discrete $(n-2)$-surface: by Proposition~\ref{propo:border.neighborhood}, since $X$ is a simplicial complex and an $n$-PCM, we obtain that $\vert\border\TC_X(h)\vert$ is equal to $\vert\TC_{\border X}(h)\vert$ and thus is a discrete $(n-2)$-surface. So, $\vert\TC_X(h)\vert$ is a smooth $(n-1)$-PCM. Fourth, we want to prove that $\border X$ is a discrete $(n-1)$-surface or a separated union of discrete $(n-1)$-surfaces. We can remark that it is equivalent to prove that, for any $h \in \border X$, $\vert \TC_{\border X}(h)\vert$ is a $(n-2)$-surface, which corresponds to the constraint $(C)$.

\medskip

The proof is done. \qed

\section{Conclusion}
\label{sec.conclusion}

We will have seen that in the specific context of simplicial complexes, discrete surfaces are equivalent to normal pseudomanifolds without border and that $n$-PCM's are equivalent to normal pseudomanifold with border. We will have also seen that it means a very strong property of these two combinatorial structures: PCM's and discrete surfaces are together the generalization of normal pseudomanifolds from simplicial complexes to posets.

\medskip

Furthermore, we will have seen a counter-example showing that there exist $n$-PCM's which are not smooth $n$-PCM's, even if they are simplicial complexes (and then supplied with many additional properties). This clearly shows how much it is important for normal pseudomanifolds to be smooth $n$-PCM's and not only $n$-PCM's (assuming that when their border is not empty).

\section{Future works}
\label{sec.future}

The theoretical results we obtain here have a direct consequence from a \emph{computational} point of view: starting from a simplicial complex $X$, and assuming that the border $\triangle X$ (empty or not) is already computed or given, checking with some algorithm whether $X$ is a normal $n$-pseudomanifold will indicate if $X$ is an $n$-PCM (when $\border X \neq \emptyset$) and if $X$ is a discrete $n$-surface (when $\border X = \emptyset$). In other words, we will be able to determine the nature of $X$ without the computationally expensive \emph{recursive} procedures induced by the definitions of PCM's or discrete surfaces. We plan then in the future to implement algorithms able to compute a border, to check whether a simplicial complex is a normal pseudomanifold, a PCM or a discrete surface, and to calculate their respective algorithmic complexities.

\bibliographystyle{plain}
\bibliography{article}

\end{document}